\input amstex
\documentstyle{amsppt}
\hsize=16cm
\vsize=23cm

\topmatter
\title
On Beilinson's Hodge and Tate conjectures for open complete intersections
\endtitle
\author{Masanori Asakura and Shuji Saito}
\endauthor

\address
{Graduate School of Mathematics,
Kyushu University 33
FUKUOKA 812-8581, JAPAN}
\endaddress
\address
{e-mail: asakura\@math.kyushu-u.ac.jp}
\endaddress
\address
{Graduate School of Mathematics,
Nagoya University, Chikusa-ku, NAGOYA, 464-8602, JAPAN}
\endaddress
\address
{e-mail: sshuji\@msb.biglobe.ne.jp}
\endaddress

\endtopmatter
\document
\NoBlackBoxes

\input amstex
\documentstyle{amsppt}
\hsize=16cm
\vsize=23cm
\document

\head \bf Contents
\endhead

\vskip 10pt\noindent
\roster
\item"\S0"
Introduction
\item"\S1"
Jacobian rings for open complete intersections
\item"\S2"
Hodge theoretic implication of generalized Jacobian rings
\item"\S3"
Beilinson's Hodge conjecture
\item"\S4"
Beilinson's Tate conjecture
\item"\S5"
Implication on injectivity of Chern class map for $K_1$ of surfaces
\item"\hbox{ }"
References
\endroster

\vskip 10pt

\input amstex
\documentstyle{amsppt}
\hsize=16cm
\vsize=23cm

\def\Th#1.{\vskip 6pt \medbreak\noindent{\bf Theorem(#1).}}
\def\Cor#1.{\vskip 6pt \medbreak\noindent{\bf Cororally(#1).}}
\def\Conj#1.{\vskip 6pt \medbreak\noindent{\bf Conjecture(#1).}}
\def\Pr#1.{\vskip 6pt \medbreak\noindent{\bf Proposition(#1).}}
\def\Lem#1.{\vskip 6pt \medbreak\noindent{\bf Lemma(#1).}}
\def\Rem#1.{\vskip 6pt \medbreak\noindent{\it Remark(#1).}}
\def\Fact#1.{\vskip 6pt \medbreak\noindent{\it Fact(#1).}}
\def\Claim#1.{\vskip 6pt \medbreak\noindent{\it Claim(#1).}}
\def\Def#1.{\vskip 6pt \medbreak\noindent{\bf Definition\bf(#1)\rm.}}

\def\qwith{\quad\hbox{with }}
\def\mathrm#1{\rm#1}

\def\isom{@>\cong>>}
\def\Spec{{\operatorname{Spec}}}

\def\dim{{\operatorname{dim}}}

\def\Coker{{\text{\rm Coker}}}
\def\dim{\hbox{\rm dim}}

\def\Im{\hbox{\rm Im}}
\def\Ker{\hbox{\rm Ker}}
\def\Coker{\hbox{\rm Coker}}
\def\min{\hbox{\rm min}}

\def\Gal{\hbox{\mathrm{Gal}}}
\def\GL{\hbox{\mathrm{GL}}}

\def\P{{\Bbb{P}}}
\def\bP{{\Bbb{P}}}

\def\cHom{{\Cal{H}}om}

\def\cF{{\Cal{F}}}
\def\cO{{\Cal{O}}}
\def\cX{{\Cal{X}}}
\def\cXB{{\Cal{X}_B}}

\def\cZ{{\Cal{Z}}}
\def\cZB{{\Cal{Z}_B}}
\def\cZBj{{\Cal{Z}_{B,j}}}

\def\wtd#1{\widetilde{#1}}

\def\l{\ell}
\def\d{{\bold{d}}}
\def\e{{\bold{e}}}

\def\ra{\rightarrow}

\def\ot{\otimes}


\def\us#1#2{\underset{#1}\to{#2}}
\def\os#1#2{\overset{#1}\to{#2}}

\def\qaq{\quad\hbox{ and }\quad}
\def\qfor{\quad\hbox{ for }}

\def\spa{\hbox{ }}
\def\scs{\spa : \spa}

\def\sp{\hbox{}}

\def\ul#1{\underline{#1}}

\def\ba{\bold{a}}
\def\bb{\bold{b}}
\def\ua{\underline{a}}
\def\ub{\underline{b}}

\def\ue{\underline{e}}
\def\ud{\underline{d}}

\def\md{\delta_{\text{min}}}

\def\WS#1{\Omega_S^{#1}}

\def\otO{\otimes_{\cO}}
\def\cU{{\Cal U}}

\def\onab{\overline{\nabla}}

\def\HO{H_{\cO}}
\def\HQ{H_{\Bbb Q}}
\def\HC{H_{\Bbb C}}
\def\HQcU#1#2{\HQ^{#1}(\cU/S)(#2)}
\def\HCcU#1{\HC^{#1}(\cU/S)}
\def\HOcU#1{\HO^{#1}(\cU/S)}

\def\HcU#1#2{H^{#1,#2}(\cU/S)}

\def\HcU#1#2{H^{#1,#2}(\cU/S)}
\def\HcXp#1#2{H^{#1,#2}(\cX/S)_{prim}}

\def\HcUan#1#2{H^{#1,#2}(\cU/S)^{an}}

\def\HQcU#1#2{\HQ^{#1}(\cU/S)(#2)}

\def\HQcU#1#2{\HQ^{#1}(\cU/S)(#2)}

\def\HCcU#1{\HC^{#1}(\cU/S)}
\def\HCcX#1{\HC^{#1}(\cX/S)}
\def\HCcU#1{\HC^{#1}(\cU/S)}
\def\HCcXp#1{\HC^{#1}(\cX/S)_{prim}}

\def\HOcU#1{\HO^{#1}(\cU/S)}
\def\HOcX#1{\HO^{#1}(\cX/S)}
\def\HOcU#1{\HO^{#1}(\cU/S)}

\def\WSan#1{\Omega_{\San}^{#1}}

\def\Xx{X_x}
\def\Zx{Z_x}
\def\Ux{U_x}
\def\xb{\overline{x}}
\def\Uxb{U_{\xb}}

\def\Zxx{Z_{x}}
\def\Zst{Z}
\def\cZst{\cZ}

\def\Yst{Y}

\def\rocXZ{\kappa_{(\cX,\cZ)}}

\def\roolog{\kappa_0^{log}}

\def\roxlog{\kappa_x^{log}}

\def\TS{\Theta_S}
\def\TxS{T_x S}

\def\UC{U_{\Bbb C}}

\def\SC{S_{\Bbb C}}

\def\WS#1{\Omega_S^{#1}}
\def\WX#1{\Omega_X^{#1}}

\def\TX{T_X}

\def\TXZx{T_{\Xx}(-\log \Zxx)}

\def\WXZ#1{\Omega_X^{#1}(\log \Zst)}

\def\WcXZ#1{\Omega_{\cX/S}^{#1}(\log \cZst)}

\def\WcXkZ#1{\Omega_{\cX/k}^{#1}(\log \cZst)}

\def\WPnY#1{\Omega_{\Bbb P^n}^{#1}(\log \Yst)}
\def\TXZ{T_X(-\log \Zst)}
\def\TcXZS{T_{\cX/S}(-\log \cZst)}

\def\fXZ#1{\phi_{X,Z}^{#1}}
\def\cXZ{\psi_{(X,Z)}}

\def\cXZx{\psi_{(\Xx,\Zx)}}

\def\eXZ{\eta_{(X,Z)}}

\def\ccXZS{c_S(\cX,\cZ)}
\def\ccXZtS{c_{\wtd{S}}(\wtd{\cX},\wtd{\cZ})}

\def\tS{\wtd{S}}

\def\scs{\hbox{ }:\hbox{ }}
\def\onab{\overline{\nabla}}

\def\dlog#1{\frac{d#1}{#1}}
\def\chB#1#2#3{ch_{B,#3}^{#1,#2}}
\def\chD#1#2#3{ch_{D,#3}^{#1,#2}}
\def\chet#1#2#3{ch_{et,#3}^{#1,#2}}
\def\chcont#1#2#3{ch_{cont,#3}^{#1,#2}}
\def\Ql{\Bbb Q_{\ell}}
\def\Zl{\Bbb Z_{\ell}}
\def\etab{\overline{\eta}}

\def\Pd{\overset{\vee}\to{\Bbb P}}
\def\kb{\overline{k}}
\def\xb{\overline{x}}

\def\pitop{\pi_1^{top}}
\def\pialg{\pi_1^{alg}}

\def\wcU{\omega_{\cU/S}(\sigma)}
\def\dcU{\delta_{\cU/S}(\sigma)}

\def\ccU{c_{\cU/S}(\sigma)}
\def\cUx{c_{\Ux}(\sigma)}
\def\cUo{c_{U}(\sigma)}
\def\cUetab{c_{U_{\etab}}(\sigma)}

\def\Dqi{\Delta_{\gamma_i}}
\def\Dq{\Delta_{\gamma}}
\def\spa{\hbox{ }}

\def\nabmo{\overline{\nabla}^{m,0}}
\def\nab#1#2{\overline{\nabla}^{#1,#2}}

\def\Pol{P}
\def\h#1#2{h_{#1}(#2)}

\def\San{S_{an}}
\def\Szar{S_{zar}}

\head \S0. Introduction. \endhead
\vskip 8pt

In his lectures in [G1], M. Green gives a lucid explanation how fruitful the
infinitesimal method in Hodge theory is in various aspects of algebraic
geometry. A significant idea is to use Koszul cohomology for Hodge-theoretic
computations. The idea originates from Griffiths work [Gri] where the
Poincar\'e residue representation of the cohomology of a hypersurface played
a crucial role in proving the infinitesimal Torelli theorem for hypersurfaces.
Since then many important applications of the idea have been made in different
geometric problems such as the generic Torelli problem and the
Noether-Lefschetz theorem and the study of algebraic cycles
(see [G1, Lectures 7 and 8]).
\vskip 3pt

In this paper we introduce \it Jacobian rings of open complete
intersections \rm and to apply it to the Beilinson's Hodge and Tate conjectures.Here, by ``open complete intersection" we mean a pair
$(X,Z=\underset{1\leq j\leq s}\to{\cup} Z_j)$ where $X$ is a smooth complete
intersection in $\Bbb P^n$ and $Z_j\subset X$ is a smooth hypersurface section
such that $Z$ is a simple normal crossing divisor on $X$.
Our Jacobian rings give an algebraic description of the infinitesimal part of
mixed Hodge structure on the cohomology $H^m(X\setminus Z,\Bbb Q)$ with
$m=\dim(X)$. It is a natural generalization of the Poincar\'e residue
representation of the cohomology of a hypersurface in [Gri].
\vskip 3pt

The Beilinson's Hodge and Tate conjectures (cf. [J1, Conjecture 8.5 and 8.6])
concern the surjectivity of the Chern class maps for open varieties.
To be more precise we let $U$ be a smooth variety over $k$.
\roster
\item
"(0-1)" (Hodge version)
When $k=\Bbb C$ the conjecture predicts the surjectivity of
the Chern class map from the higher Chow group to the Betti cohomology
(cf. [Bl] and [Sch])
$$ \chB i j U \scs CH^j(U,2j-i)\otimes\Bbb Q \to (2\pi\sqrt{-1})^j
W_{2j}H_B^i(U(\Bbb C),\Bbb Q)\cap F^jH_B^i(U(\Bbb C),\Bbb C)$$
where $W_*$ (resp. $F^*$) denotes the weight (resp. Hodge) filtration of
the mixed Hodge structure on the Betti cohomology group.
\item
"(0-2)" (Tate version)
When $k$ is a finite extension of $\Bbb Q$ the conjecture predicts the
surjectivity of the Chern class map from the higher Chow group to the etale
cohomology
$$ \chet i j U \scs CH^j(U,2j-i)\otimes\Ql \to
H^i_{et}(U\times_k \kb,\Ql(j))^{Gal(\kb/k)}$$
where $\kb$ is an algebraic closure of $k$.
\endroster

\Conj 0-1. (A.Beilinson [Bei]) \it The above maps are surjective in case $i=j$.
\rm\vskip 5pt

The following are some remarks on the conjecture.
\roster
\item
"$(i)$"
In case that $U$ is proper and $i=2j$,
the surjectivity of the above maps is equivalent to the Hodge and Tate
conjectures for algebraic cycles of codimension $j$ on $U$.
Thus Conj.(2-1) is analous to the Hodge and Tate conjectures for
algebraic cycles on a projective smooth variety.
\item
"$(ii)$"
Let $X$ be a projective smooth curve and $U\subset X$ be a non-empty open
subset. Then the surjectivity of $\chB 1 1 U$ (resp. $\chet 1 1 U$) for
follows from the Abel's theorem. Indeed it is equivalent to the injectivity
(modulo torsion) of the Abel-Jacobi map (resp. $\ell$-adic Abel-Jacobi map)
$$ CH_0(X)^{deg=0} \to Jac(X)
\quad (\text{resp. }
CH_0(X)^{deg=0} \to H_{gal}^1(k,H^1(X\times_k\kb,\Ql(1)))$$
restricted on the subspace of $CH_0(X)^{deg=0}$ generated by cycles supported
on $X\setminus U$.
Here $CH_0(X)^{deg=0}$ is the group of zero-cycles on $X$ of degree zero
modulo rational equivalence and $Jac(X)$ is the Jacobian variety of $X$.
\item
"$(iii)$"
When $X$ is a projective smooth surface and $U\subset X$ is the complement
of a simple normal crossing divisor $Z\subset X$, then the surjectivity of
$\chB 22 U$ or $\chet 22 U$ has an implication on the injectivity of the Chern
class map for $ CH^2(X,1)$, which can be viewed as an analogue of the Abel's
theorem for $K_1$ of surfaces. The detail will be discussed in \S5.
\item
"$(iv)$"
As a naive generalization of the Beilinson's conjecture, one may ask if
the Chern class maps in (0-1) and (0-2) are surjective for any $i,j\geq 0$.
Jannsen ([J1, 9.11]) has shown that the map in (0-1) in case $i=1$,
$j=2q-1\geq3$ is not surjective in general by using a theorem of Mumford [Mu],
which implies the Abel-Jacobi for cycles of codimension$\geq 2$ is not
injective even modulo torsion.
\endroster
\vskip 5pt

In order to state the main result on the Beilinson conjectures
we fix a field $k$ of characteristic zero and
a non-singular quasi-projective variety $S$ over $k$. We also fix integers
$d_1,\dots,d_r,e_1,\dots,e_s\geq 1$.
Assume that we are given schemes over $S$
$$\Bbb P^n_S \hookleftarrow \cX \hookleftarrow
\cZ=\underset{1\leq j\leq s}\to{\cup} \cZ_j$$
whose fibers are open complete intersections.
We assume that the fibers of $\cX/S$ are smooth complete intersection of
multi-degree $(d_1,\dots,d_r)$ in $\Bbb P^n$ and that those of
$\cZ_j\subset\cX$ are smooth hypersurface section of degree $e_j$.
Let $f:\cX\to S$ be the natural morphism and write $\cU=\cX \setminus \cZ$.
Let $\Ux$ denote the fiber of $\cU=\cX\setminus\cZ$ over $x\in S$.
In \S2 we will introduce an invariant $\ccXZS$ that measures
the ``generality" of the family (0-1), or
how many independent parameters $S$ contains (cf. Rem.(2-2)).

\Th 0-1. \it Assume $\underset{1\leq i\leq r}\to{\sum} d_i \geq n+1+\ccXZS$.
Let $m=n-r\geq 1$.
\roster
\item
Assume $k=\Bbb C$. There exists $E\subset S(\Bbb C)$ which is
the union of countable many proper analytic subset of $S(\Bbb C)$ such that
$\chB m m \Ux$ is surjective for $\forall x\in S(\Bbb C)\setminus E$.
\item
Assume that $k$ is a finite extension of $\Bbb Q$ and
$S(k)\not=\emptyset$.
Let $\pi: S \to \Bbb P_k^N$ be a dominant quasi-finite morphism.
There exist a subset $H\subset\Bbb P_k^N(k)$ such that:
\vskip 4pt\noindent
$(i)$
$\chet m m \Ux$ is surjective for any closed point $x\in S$ such that
$\pi(x)\in H$.
\vskip 4pt\noindent
$(ii)$
Let $\Sigma$ be any finite set of primes of $k$ and let $k_v$ be the
completion of $k$ at $v\in \Sigma$. Then the image of $H$ in
$\prod_{v\in \Sigma} \Bbb P_k^N(k_v)$ is dense.
\endroster
\rm
\vskip 6pt

Indeed the target spaces of the maps $\chB m m {\Ux}$ and
$\chet m m {\Ux}$ are non-zero if $s\geq n-r+1$ (recall that $s$ is the
number of the irreducible components of $\cZ$) and we give explicit
elements in $CH^m(\Ux,m)$ whose images span them.
\vskip 6pt

Now we explain how the paper is organized.
In \S1 we state the fundamental results on the generalized Jacobian rings,
the duality theorem and the symmetrizer lemma.
The proof is given in another paper [AS1]. It is based on the basic techniques
to compute Koszul cohomology developed by M. Green ([G2] and [G3]).
In \S2 we give a Hodge theoretic implication of the results in \S1 which
plays a crucial role in the proof of Th.(0-1).
Th.(0-1)(1) is proven in \S3 and Th.(0-1)(2) is proven in \S4
by using the results in \S2.
In \S5 we explain an implication of the Beilinson's conjectures on the
injectivity of Chern class maps for $K_1$ of surfaces.

\vskip 20pt

\input amstex
\documentstyle{amsppt}
\hsize=16cm
\vsize=23cm

\def\Th#1.{\vskip 6pt \medbreak\noindent{\bf Theorem(#1).}}
\def\Cor#1.{\vskip 6pt \medbreak\noindent{\bf Cororally(#1).}}
\def\Conj#1.{\vskip 6pt \medbreak\noindent{\bf Conjecture(#1).}}
\def\Pr#1.{\vskip 6pt \medbreak\noindent{\bf Proposition(#1).}}
\def\Lem#1.{\vskip 6pt \medbreak\noindent{\bf Lemma(#1).}}
\def\Rem#1.{\vskip 6pt \medbreak\noindent{\it Remark(#1).}}
\def\Fact#1.{\vskip 6pt \medbreak\noindent{\it Fact(#1).}}
\def\Claim#1.{\vskip 6pt \medbreak\noindent{\it Claim(#1).}}
\def\Def#1.{\vskip 6pt \medbreak\noindent{\bf Definition\bf(#1)\rm.}}

\def\qwith{\quad\hbox{with }}
\def\mathrm#1{\rm#1}

\def\isom{@>\cong>>}
\def\Spec{{\operatorname{Spec}}}

\def\dim{{\operatorname{dim}}}

\def\Coker{{\text{\rm Coker}}}
\def\dim{\hbox{\rm dim}}

\def\Im{\hbox{\rm Im}}
\def\Ker{\hbox{\rm Ker}}
\def\Coker{\hbox{\rm Coker}}
\def\min{\hbox{\rm min}}

\def\Gal{\hbox{\mathrm{Gal}}}
\def\GL{\hbox{\mathrm{GL}}}

\def\P{{\Bbb{P}}}
\def\bP{{\Bbb{P}}}

\def\cHom{{\Cal{H}}om}

\def\cF{{\Cal{F}}}
\def\cO{{\Cal{O}}}
\def\cX{{\Cal{X}}}
\def\cXB{{\Cal{X}_B}}

\def\cZ{{\Cal{Z}}}
\def\cZB{{\Cal{Z}_B}}
\def\cZBj{{\Cal{Z}_{B,j}}}

\def\wtd#1{\widetilde{#1}}

\def\l{\ell}
\def\d{{\bold{d}}}
\def\e{{\bold{e}}}

\def\ra{\rightarrow}

\def\ot{\otimes}


\def\us#1#2{\underset{#1}\to{#2}}
\def\os#1#2{\overset{#1}\to{#2}}

\def\qaq{\quad\hbox{ and }\quad}
\def\qfor{\quad\hbox{ for }}

\def\spa{\hbox{ }}
\def\scs{\spa : \spa}

\def\sp{\hbox{}}

\def\ul#1{\underline{#1}}

\def\ba{\bold{a}}
\def\bb{\bold{b}}
\def\ua{\underline{a}}
\def\ub{\underline{b}}

\def\ue{\underline{e}}
\def\ud{\underline{d}}

\def\md{\delta_{\text{min}}}

\def\WS#1{\Omega_S^{#1}}

\def\otO{\otimes_{\cO}}
\def\cU{{\Cal U}}

\def\onab{\overline{\nabla}}

\def\HO{H_{\cO}}
\def\HQ{H_{\Bbb Q}}
\def\HC{H_{\Bbb C}}
\def\HQcU#1#2{\HQ^{#1}(\cU/S)(#2)}
\def\HCcU#1{\HC^{#1}(\cU/S)}
\def\HOcU#1{\HO^{#1}(\cU/S)}

\def\HcU#1#2{H^{#1,#2}(\cU/S)}

\def\HcU#1#2{H^{#1,#2}(\cU/S)}
\def\HcXp#1#2{H^{#1,#2}(\cX/S)_{prim}}

\def\HcUan#1#2{H^{#1,#2}(\cU/S)^{an}}

\def\HQcU#1#2{\HQ^{#1}(\cU/S)(#2)}

\def\HQcU#1#2{\HQ^{#1}(\cU/S)(#2)}

\def\HCcU#1{\HC^{#1}(\cU/S)}
\def\HCcX#1{\HC^{#1}(\cX/S)}
\def\HCcU#1{\HC^{#1}(\cU/S)}
\def\HCcXp#1{\HC^{#1}(\cX/S)_{prim}}

\def\HOcU#1{\HO^{#1}(\cU/S)}
\def\HOcX#1{\HO^{#1}(\cX/S)}
\def\HOcU#1{\HO^{#1}(\cU/S)}

\def\WSan#1{\Omega_{\San}^{#1}}

\def\Xx{X_x}
\def\Zx{Z_x}
\def\Ux{U_x}
\def\xb{\overline{x}}
\def\Uxb{U_{\xb}}

\def\Zxx{Z_{x}}
\def\Zst{Z}
\def\cZst{\cZ}

\def\Yst{Y}

\def\rocXZ{\kappa_{(\cX,\cZ)}}

\def\roolog{\kappa_0^{log}}

\def\roxlog{\kappa_x^{log}}

\def\TS{\Theta_S}
\def\TxS{T_x S}

\def\UC{U_{\Bbb C}}

\def\SC{S_{\Bbb C}}

\def\WS#1{\Omega_S^{#1}}
\def\WX#1{\Omega_X^{#1}}

\def\TX{T_X}

\def\TXZx{T_{\Xx}(-\log \Zxx)}

\def\WXZ#1{\Omega_X^{#1}(\log \Zst)}

\def\WcXZ#1{\Omega_{\cX/S}^{#1}(\log \cZst)}

\def\WcXkZ#1{\Omega_{\cX/k}^{#1}(\log \cZst)}

\def\WPnY#1{\Omega_{\Bbb P^n}^{#1}(\log \Yst)}
\def\TXZ{T_X(-\log \Zst)}
\def\TcXZS{T_{\cX/S}(-\log \cZst)}

\def\fXZ#1{\phi_{X,Z}^{#1}}
\def\cXZ{\psi_{(X,Z)}}

\def\cXZx{\psi_{(\Xx,\Zx)}}

\def\eXZ{\eta_{(X,Z)}}

\def\ccXZS{c_S(\cX,\cZ)}
\def\ccXZtS{c_{\wtd{S}}(\wtd{\cX},\wtd{\cZ})}

\def\tS{\wtd{S}}

\def\scs{\hbox{ }:\hbox{ }}
\def\onab{\overline{\nabla}}

\def\dlog#1{\frac{d#1}{#1}}
\def\chB#1#2#3{ch_{B,#3}^{#1,#2}}
\def\chD#1#2#3{ch_{D,#3}^{#1,#2}}
\def\chet#1#2#3{ch_{et,#3}^{#1,#2}}
\def\chcont#1#2#3{ch_{cont,#3}^{#1,#2}}
\def\Ql{\Bbb Q_{\ell}}
\def\Zl{\Bbb Z_{\ell}}
\def\etab{\overline{\eta}}

\def\Pd{\overset{\vee}\to{\Bbb P}}
\def\kb{\overline{k}}
\def\xb{\overline{x}}

\def\pitop{\pi_1^{top}}
\def\pialg{\pi_1^{alg}}

\def\wcU{\omega_{\cU/S}(\sigma)}
\def\dcU{\delta_{\cU/S}(\sigma)}

\def\ccU{c_{\cU/S}(\sigma)}
\def\cUx{c_{\Ux}(\sigma)}
\def\cUo{c_{U}(\sigma)}
\def\cUetab{c_{U_{\etab}}(\sigma)}

\def\Dqi{\Delta_{\gamma_i}}
\def\Dq{\Delta_{\gamma}}
\def\spa{\hbox{ }}

\def\nabmo{\overline{\nabla}^{m,0}}
\def\nab#1#2{\overline{\nabla}^{#1,#2}}

\def\Pol{P}
\def\h#1#2{h_{#1}(#2)}

\def\San{S_{an}}
\def\Szar{S_{zar}}

\head \S1. Jacobian rings for open complete intersections.
\endhead
\vskip 8pt

The purpose of this section is to introduce Jacobian rings for
open compolete intersections and state their fundamental properties.
Throughout the whole paper, we fix integers $r,s \geq 0$ with $r+s\geq 1$,
$n\geq 2$ and $d_1, \cdots, d_r$, $e_1, \cdots, e_s \geq 1$. We put
$$
\d=\sum_{i=1}^{r}d_i, \quad \e=\sum_{j=1}^{s}e_j,\quad
\md=\underset{\underset{1\leq j\leq s}\to{1\leq i\leq r}}\to{\min}\{d_i,e_j\},
\quad d_{max}=\underset{1\leq i\leq r}\to{\max}\{d_i\},
\quad e_{max}=\underset{1\leq j\leq s}\to{\max}\{e_j\}.
$$
We also fix a field $k$ of characteristic zero.
Let $\Pol=k[X_0,\dots,X_n]$ be the polynomial ring over $k$ in $n+1$
variables.
Denote by $\Pol^l\subset \Pol$ the subspace of the homogeneous
polynomials of degree $l$.
Let $A$ be a polynomial ring over $\Pol$ with indeterminants
$\mu_1, \cdots, {\mu}_r$, $\lambda_1,\cdots, \lambda_s$.
We use the multi-index notation
$$\mu^{\underline{a}}=\mu_1^{a_1}\cdots\mu_r^{a_r} \qaq
 \lambda^{\underline{b}}=\lambda_1^{b_1}\cdots\lambda_s^{b_s} \qfor
\underline{a}=(a_1,\cdots,a_r) \in {\Bbb Z}^{\oplus r}_{\geq 0},\spa
\underline{b}=(b_1,\cdots,b_s) \in {\Bbb Z}^{\oplus s}_{\geq 0}.$$
For $q \in {\Bbb Z}$ and $\ell \in {\Bbb Z}$, we write
$$
 A_q(\ell)=
\us{\ba+\bb=q}{\oplus}
\Pol^{\ua\ud+\ub\ue+\ell} \cdot
\mu^{\underline{a}}\lambda^{\underline{b}}\quad
(\ba=\sum_{i=1}^r a_i,\sp
\bb=\sum_{j=1}^s b_j,\sp
\ua\ud=\sum_{i=1}^r a_id_i,\sp
\ub\ue=\sum_{j=1}^s b_je_j)
$$
By convention $A_q(\ell)=0$ if $q<0$.

\Def1-1.
For $\ul{F}=(F_1, \cdots, F_r)$, $\ul{G}=(G_1, \cdots, G_s)$ with
$F_i \in \Pol^{d_i}$, $G_j \in \Pol^{e_j}$, we define the {\bf Jacobian ideal}
$J(\ul{F},\ul{G})$ to be the ideal of $A$ generated by
$$\sum_{1\leq i\leq r}\frac{\partial F_i}{\partial X_k}\mu_i+
\sum_{1\leq j\leq s}\frac{\partial G_j}{\partial X_k}\lambda_j,
\quad F_i,\quad G_j\lambda_j \quad
(1 \leq i \leq r,\sp 1 \leq j \leq s,\sp 0 \leq k \leq n).$$
The quotient ring $B=B(\ul{F},\ul{G})=A/J(\ul{F},\ul{G})$ is called
the {\bf Jacobian ring of} ($\ul{F},\ul{G}$). We denote
$$
B_q(\ell) =B_q(\ell)(\ul{F},\ul{G})=A_q(\ell)/A_q(\ell) \cap J(\ul{F},\ul{G}).
$$
\vskip 6pt

\Def 1-2. Suppose $n\geq r+1$.
Let $\Bbb{P}^n=\text{\rm Proj }\Pol$ be the projective space over $k$.
Let $X\subset \P^n$ be defined by $F_1=\cdots=F_r=0$ and let
$Z_j\subset X$ be defined by $G_j=F_1=\cdots=F_r=0$ for $1 \leq j \leq s$.
We also call $B(\ul{F},\ul{G})$ the Jacobian ring of the pair
$(X,\Zst={\cup}_{1\leq j\leq s} Z_j)$ and denote
$B(\ul{F},\ul{G})=B(X,\Zst)$ and $J(\ul{F},\ul{G})=J(X,\Zst)$.

\vskip 6pt
In what follows we fix $\ul{F}$ and $\ul{G}$ as Def.(1-1) and assume
the condition
$$ \text{$F_i=0$ ($1\leq i\leq r$) and
$G_j=0$ ($1\leq j\leq s$) intersect transversally in $\Bbb{P}^n$.}
\leqno(1-1)$$
We mension three main theorems.
The first main theorem concerns with the geometric meaning of Jacobian rings.

\Th I. \it
Suppose $n\geq r+1$. Let $X$ and $Z$ be as Definition (1-2).
\vskip 4pt\noindent
(1)
For intergers $0 \leq q\leq n-r$ and $\ell\geq 0$ there is a natural
isomorphism
$$\fXZ q \scs B_q(\d+\e-n-1+\ell) \isom H^{q}(X, \WXZ {n-r-q}(\ell))_{prim}.$$
Here $\WXZ p$ is the sheaf of algebraic differential $q$-forms on $X$
with logarithmic poles along $\Zst$ and `$prim$' means the primitive part:
$$ H^{q}(X, \WXZ {p}(\ell))_{prim}=
\left\{\aligned
& \Coker(H^q(\bP^n,\Omega_{\bP^n}^q) \to H^{q}(X, \WX {q}))
\quad\text{ if $q=p$ and $s=\ell=0$,}\\
& H^{q}(X, \WXZ {p}(\ell))
\quad\text{ otherwise.}\\
\endaligned\right.$$
\vskip 4pt\noindent
(2)
There is a natural map
$$\cXZ : B_1(0) \longrightarrow H^1(X, \TXZ)_{alg}\subset H^1(X, \TXZ) $$
which is an isomorphism if $\dim(X)\geq 2$.
Here $\TXZ$ is the $\cO_X$-dual of $\WXZ 1$ and the group on the right hand
side is defined in Def.(1-3) below. The following map
$$
H^1(X, \TXZ) \otimes H^{q}(X, \WXZ p)
\longrightarrow H^{q+1}(X, \WXZ {p-1}).
$$
induced by the cup-product and the contraction
$\TXZ\ot\WXZ p\to\WXZ {p-1}$
is compatible through $\cXZ$ with the ring multiplication up to scalar.
\rm\vskip 6pt

Roughly speaking, the generalized Jacobian rings describe the infinitesimal
part of the Hodge structures of open
variety $X \setminus Z$,
and the cup-product with Kodaira-Spencer class
coincides with the ring multiplication up to non-zero scalar.
This result was originally invented by P. Griffiths in case of hypersurfaces
and generalized to complete intersections by Konno [K].
Our result is a further generalization.

\Def 1-3. \it Let the assumption be as in Th.(I).
We define $H^1(X, \TXZ)_{alg}$ to be
the kernel of the composite map
$$H^1(X,\TXZ)\to H^1(X,\TX) \to H^2(X,\cO_X),$$
where the second map is induced by the cup product with the class
$c_1(\cO_{X}(1))\in H^1(X,\WX 1)$ and the contraction
$\TX\otimes \WX 1 \to \cO_X$.
It can be seen that
$$ \dim_k(H^1(X, \TXZ)/H^1(X, \TXZ)_{alg})=\left\{
\aligned
1 & \text{ if $X$ is a $K3$ surface,}\\
0 & \text{ otherwise.}\\
\endaligned \right.
$$

\vskip 6pt\rm

The second main theorem is the duality theorem for the generalized
Jacobian rings.

\Th II. \it
(1)
There is a natural map (called the trace map)
$$\tau\sp:\sp B_{n-r}(2(\d-n-1)+\e) \to k.$$
Let
$$ \h p \ell\sp:\sp B_p(\d-n-1+\l) \ra B_{n-r-p}(\d+\e-n-1-\l)^*$$
be the map induced by the following pairing induced by
the multiplication
$$B_{p}(\d-n-1+\l) \otimes B_{n-r-p}(\d+\e-n-1-\l)\to  B_{n-r}(2(\d-n-1)+\e)
@>\tau>> k.$$
When $r>n$ we define $ \h p \ell$ to be the zero map by convention.
\vskip 4pt\noindent
(2)
The map $\h p \ell$ is an isomorphism in either of the following cases.
\roster
\item"$(i)$"  $s\geq 1$ and $p<n-r$ and $\ell<e_{max}$.
\item"$(ii)$"  $s\geq 1$ and $0\leq \ell\leq e_{max}$ and $r+s\leq n$.
\item"$(iii)$" $s=\ell=0$ and either $n-r\geq 1$ or $n-r=p=0$.
\endroster
\vskip 4pt\noindent
(3)
The map $\h {n-r} \ell$ is injective if $s\geq 1$ and $\ell<e_{max}$.
\vskip 6pt\rm

We have the following auxiliary result on the duality.

\Th II'. \it
Assume $n-r\geq 1$ and consider the composite map
$$ \eXZ \scs H^{0}(X, \WXZ {n-r})@>{(\fXZ 0)^{-1}}>>
B_0(\d+\e-n-1) @>{{\h {n-r} 0}^*}>>  B_{n-r}(\d-n-1)^*$$
where the second map is the dual of $\h {n-r} 0$.
Then $\eXZ$ is surjective and we have (cf. Def.(1-4) below)
$$ \Ker(\eXZ)=\wedge_X^{n-r}(G_1,\dots,G_s).$$
\rm\vskip 6pt

\Def 1-4. \it
Let $G_1,\dots,G_s$ be as in Def.(1-1) and let $Y_j\subset \P^n$ be
the smooth hypersurface defined by $G_j=0$.
Let $X\subset \P^n$ be a smooth projective variety such that
$Y_j$ ($1\leq j\leq s$) and $X$ intersect transversally.
Put $Z_j=X\cap Y_j$.
Take an integer $q$ with $0\leq q \leq s-1$.
For integers $1\leq j_1<\cdots< j_{q+1}\leq s$, let
$$ \omega_X(j_1,\dots,j_{q+1})\in H^{0}(X, \WXZ {q})
\quad (\Zst=\underset{1\leq j\leq s}\to{\Sigma} Z_j)$$
be the restriction of
$$ \sum_{\nu=1}^{q+1}(-1)^{\nu-1} e_{j_{\nu}}
\frac{dG_{j_1}}{G_{j_1}}\wedge\cdots\wedge\widehat
{\frac{dG_{j_\nu}}{G_{j_\nu}}}
\wedge\cdots\wedge{\frac{dG_{j_{q+1}}}{G_{j_{q+1}}}} \sp\in
H^0(\P^n,\WPnY q)$$
where
$\Yst=\underset{1\leq j\leq s}\to{\Sigma} Y_j\subset \P^n$.
We let
$$ \wedge_X^q(G_1,\dots,G_s)\subset H^{0}(X, \WXZ {q})$$
be the subspace generated by $\omega_X(j_1,\dots,j_{q+1})$.
For $1\leq j_1<\cdots< j_{q}\leq s-1$ we have
$$  e_s\cdot \omega_X(j_1,\dots,j_{q},s)=
\frac{dg_{j_1}}{g_{j_1}}\wedge\cdots\wedge
{\frac{dg_{j_{q}}}{g_{j_{q}}}}
\qwith g_j=(G_j^{e_s}/G_s^{e_j})_{|X}\in \Gamma(U,\cO_{U}^*)
\quad (U=X\setminus Z)$$
and $\omega_X(j_1,\dots,j_{q},s)$ with $1\leq j_1<\cdots< j_{q}\leq s-1$
form a basis of $\wedge_X^q(G_1,\dots,G_s)$.
\rm\vskip 6pt

Our last main theorem is the generalization of Donagi's symmetrizer lemma [Do]
(see also [DG], [Na] and [N]) to the case of open complete intersections at
higher degrees.

\Th III. \it Assume $s\geq 1$.
Let $V \subset B_1(0)$ is a subspace of codimension $c\geq 0$.
Then the Koszul complex
$$
B_p(\l) \otimes \os{q+1}{\wedge}V \rightarrow
B_{p+1}(\l) \otimes \os{q}{\wedge}V \rightarrow
B_{p+2}(\l) \otimes \os{q-1}{\wedge}V
$$
is exact if one of the following conditions is  satisfied.
\roster
\item"$(i)$"
$p\geq 0$, $q=0$ and $\md p+\ell\geq c$.
\item"$(ii)$"
$p\geq 0$, $q=1$ and $\md p+\ell\geq 1+c$ and $\md(p+1)+\ell\geq d_{max}+c$.
\item"$(iii)$"
$p\geq 0$, $\md(r+p)+\ell\geq \d+q+c$, $\d+e_{max}-n-1>\ell\geq \d-n-1$
and either $r+s\leq n+2$ or $p\leq n-r-[q/2]$, where
$[*]$ denotes the Gaussian symbol.
\endroster
\rm\vskip 6pt

\rm

\vskip 20pt

\input amstex
\documentstyle{amsppt}
\hsize=16cm
\vsize=23cm

\def\Th#1.{\vskip 6pt \medbreak\noindent{\bf Theorem(#1).}}
\def\Cor#1.{\vskip 6pt \medbreak\noindent{\bf Cororally(#1).}}
\def\Conj#1.{\vskip 6pt \medbreak\noindent{\bf Conjecture(#1).}}
\def\Pr#1.{\vskip 6pt \medbreak\noindent{\bf Proposition(#1).}}
\def\Lem#1.{\vskip 6pt \medbreak\noindent{\bf Lemma(#1).}}
\def\Rem#1.{\vskip 6pt \medbreak\noindent{\it Remark(#1).}}
\def\Fact#1.{\vskip 6pt \medbreak\noindent{\it Fact(#1).}}
\def\Claim#1.{\vskip 6pt \medbreak\noindent{\it Claim(#1).}}
\def\Def#1.{\vskip 6pt \medbreak\noindent{\bf Definition\bf(#1)\rm.}}

\def\qwith{\quad\hbox{with }}
\def\mathrm#1{\rm#1}

\def\isom{@>\cong>>}
\def\Spec{{\operatorname{Spec}}}

\def\dim{{\operatorname{dim}}}

\def\Coker{{\text{\rm Coker}}}
\def\dim{\hbox{\rm dim}}

\def\Im{\hbox{\rm Im}}
\def\Ker{\hbox{\rm Ker}}
\def\Coker{\hbox{\rm Coker}}
\def\min{\hbox{\rm min}}

\def\Gal{\hbox{\mathrm{Gal}}}
\def\GL{\hbox{\mathrm{GL}}}

\def\P{{\Bbb{P}}}
\def\bP{{\Bbb{P}}}

\def\cHom{{\Cal{H}}om}

\def\cF{{\Cal{F}}}
\def\cO{{\Cal{O}}}
\def\cX{{\Cal{X}}}
\def\cXB{{\Cal{X}_B}}

\def\cZ{{\Cal{Z}}}
\def\cZB{{\Cal{Z}_B}}
\def\cZBj{{\Cal{Z}_{B,j}}}

\def\wtd#1{\widetilde{#1}}

\def\l{\ell}
\def\d{{\bold{d}}}
\def\e{{\bold{e}}}

\def\ra{\rightarrow}

\def\ot{\otimes}


\def\us#1#2{\underset{#1}\to{#2}}
\def\os#1#2{\overset{#1}\to{#2}}

\def\qaq{\quad\hbox{ and }\quad}
\def\qfor{\quad\hbox{ for }}

\def\spa{\hbox{ }}
\def\scs{\spa : \spa}

\def\sp{\hbox{}}

\def\ul#1{\underline{#1}}

\def\ba{\bold{a}}
\def\bb{\bold{b}}
\def\ua{\underline{a}}
\def\ub{\underline{b}}

\def\ue{\underline{e}}
\def\ud{\underline{d}}

\def\md{\delta_{\text{min}}}

\def\WS#1{\Omega_S^{#1}}

\def\otO{\otimes_{\cO}}
\def\cU{{\Cal U}}

\def\onab{\overline{\nabla}}

\def\HO{H_{\cO}}
\def\HQ{H_{\Bbb Q}}
\def\HC{H_{\Bbb C}}
\def\HQcU#1#2{\HQ^{#1}(\cU/S)(#2)}
\def\HCcU#1{\HC^{#1}(\cU/S)}
\def\HOcU#1{\HO^{#1}(\cU/S)}

\def\HcU#1#2{H^{#1,#2}(\cU/S)}

\def\HcU#1#2{H^{#1,#2}(\cU/S)}
\def\HcXp#1#2{H^{#1,#2}(\cX/S)_{prim}}

\def\HcUan#1#2{H^{#1,#2}(\cU/S)^{an}}

\def\HQcU#1#2{\HQ^{#1}(\cU/S)(#2)}

\def\HQcU#1#2{\HQ^{#1}(\cU/S)(#2)}

\def\HCcU#1{\HC^{#1}(\cU/S)}
\def\HCcX#1{\HC^{#1}(\cX/S)}
\def\HCcU#1{\HC^{#1}(\cU/S)}
\def\HCcXp#1{\HC^{#1}(\cX/S)_{prim}}

\def\HOcU#1{\HO^{#1}(\cU/S)}
\def\HOcX#1{\HO^{#1}(\cX/S)}
\def\HOcU#1{\HO^{#1}(\cU/S)}

\def\WSan#1{\Omega_{\San}^{#1}}

\def\Xx{X_x}
\def\Zx{Z_x}
\def\Ux{U_x}
\def\xb{\overline{x}}
\def\Uxb{U_{\xb}}

\def\Zxx{Z_{x}}
\def\Zst{Z}
\def\cZst{\cZ}

\def\Yst{Y}

\def\rocXZ{\kappa_{(\cX,\cZ)}}

\def\roolog{\kappa_0^{log}}

\def\roxlog{\kappa_x^{log}}

\def\TS{\Theta_S}
\def\TxS{T_x S}

\def\UC{U_{\Bbb C}}

\def\SC{S_{\Bbb C}}

\def\WS#1{\Omega_S^{#1}}
\def\WX#1{\Omega_X^{#1}}

\def\TX{T_X}

\def\TXZx{T_{\Xx}(-\log \Zxx)}

\def\WXZ#1{\Omega_X^{#1}(\log \Zst)}

\def\WcXZ#1{\Omega_{\cX/S}^{#1}(\log \cZst)}

\def\WcXkZ#1{\Omega_{\cX/k}^{#1}(\log \cZst)}

\def\WPnY#1{\Omega_{\Bbb P^n}^{#1}(\log \Yst)}
\def\TXZ{T_X(-\log \Zst)}
\def\TcXZS{T_{\cX/S}(-\log \cZst)}

\def\fXZ#1{\phi_{X,Z}^{#1}}
\def\cXZ{\psi_{(X,Z)}}

\def\cXZx{\psi_{(\Xx,\Zx)}}

\def\eXZ{\eta_{(X,Z)}}

\def\ccXZS{c_S(\cX,\cZ)}
\def\ccXZtS{c_{\wtd{S}}(\wtd{\cX},\wtd{\cZ})}

\def\tS{\wtd{S}}

\def\scs{\hbox{ }:\hbox{ }}
\def\onab{\overline{\nabla}}

\def\dlog#1{\frac{d#1}{#1}}
\def\chB#1#2#3{ch_{B,#3}^{#1,#2}}
\def\chD#1#2#3{ch_{D,#3}^{#1,#2}}
\def\chet#1#2#3{ch_{et,#3}^{#1,#2}}
\def\chcont#1#2#3{ch_{cont,#3}^{#1,#2}}
\def\Ql{\Bbb Q_{\ell}}
\def\Zl{\Bbb Z_{\ell}}
\def\etab{\overline{\eta}}

\def\Pd{\overset{\vee}\to{\Bbb P}}
\def\kb{\overline{k}}
\def\xb{\overline{x}}

\def\pitop{\pi_1^{top}}
\def\pialg{\pi_1^{alg}}

\def\wcU{\omega_{\cU/S}(\sigma)}
\def\dcU{\delta_{\cU/S}(\sigma)}

\def\ccU{c_{\cU/S}(\sigma)}
\def\cUx{c_{\Ux}(\sigma)}
\def\cUo{c_{U}(\sigma)}
\def\cUetab{c_{U_{\etab}}(\sigma)}

\def\Dqi{\Delta_{\gamma_i}}
\def\Dq{\Delta_{\gamma}}
\def\spa{\hbox{ }}

\def\nabmo{\overline{\nabla}^{m,0}}
\def\nab#1#2{\overline{\nabla}^{#1,#2}}

\def\Pol{P}
\def\h#1#2{h_{#1}(#2)}

\def\San{S_{an}}
\def\Szar{S_{zar}}

\head \S2. Hodge theoretic implication of generalized Jacobian rings. \endhead
\vskip 8pt

Let the assumption be as in \S1. We fix a non-singular affine algebraic
variety $S$ over $k$ and the following schemes over $S$
$$\Bbb P^n_S \hookleftarrow \cX \overset{i}\to\hookleftarrow
\cZ=\underset{1\leq j\leq s}\to{\cup} \cZ_j\leqno(2-1)$$
whose fibers are as in Def.(1-2).
Let $f:\cX\to S$ be the natural morphism and write $\cU=\cX \setminus \cZ$.
For integers $p,q$ we introduce the following sheaf on $\Szar$
$$\HcU p q=R^q f_* \WcXZ p,$$
where $\WcXZ p=\os{p}{\wedge}\WcXZ 1$ with $\WcXZ 1$, the sheaf of relative
differentials on $\cX$ over $S$ with logarithmic poles along $\cZ$.
In case $s\geq 1$ the Lefschetz theory implies $\HcU p q=0$
if $p+q\not=n-r$.
In case $s=0$ it implies $\HcXp p q=0$
if $p+q\not=n-r$ where ``$prim$"" denotes the primitive part (cf. Th.(I)(1)).
The results in \S1 implies that under an appropriate
numerical condition on $d_i$ and $e_j$ we can control the cohomology of the
following Koszul complex
$$ \WS {q-1}\ot \HcU {a+2}{b-2} @>\onab>>
 \WS {q}\ot \HcU {a+1}{b-1} @>\onab>>
 \WS {q+1}\ot \HcU {a}{b}.$$
Here $\onab$ is induced by the Kodaira-Spencer map
$$ \rocXZ\scs \TS \to R^1f_* \TcXZS,\leqno(2-2)$$
with $\TS=\cHom_{\cO_S}(\WS 1,\cO_S)$ and
$\TcXZS=\cHom_{\cO_{\cX}}(\WcXZ 1,\cO_{\cX})$, and the map
$$ R^1f_* \TcXZS \otimes R^{b-1} f_* \WcXZ {a+1} \to R^b f_* \WcXZ a$$
induced by the cup product and
$\TcXZS \otimes \WcXZ {a+1} \to \WcXZ a$, the contraction.
For the application to the Beilinsons conjectures the kernel of
the following map plays a crucial role
$$  \nab p q \scs \HcU p q @>\onab>> \WS {1}\otO \HcU {p-1}{q+1}.$$
In case $s=0$ we let $\nab p q$ denote the primitive part of the above map
$$ \HcXp p q @>\onab>> \WS {1}\otO \HcXp {p-1}{q+1}.$$
The key result is the following. The notations will be introduced later
(cf. Def.(2-1) and (2-2)).

\Th 2-1. \it Assume $p+q=m=n-r\geq 1$.
\roster
\item
Assuming $1\leq p\leq m-1$ and
$\md(p-1)+\d\geq n+1+\ccXZS$, we have $\Ker(\nab pq)=0.$
\item
Assume $\md(n-r-1)+\d\geq n+1+\ccXZS.$
Then $\Ker(\nabmo)$ is generated as $\cO_S$-module by
$\wcU$ with $\sigma\in J$. In particular $\Ker(\nabmo)=0$ if $s\leq n-r$.
\endroster
\rm
\vskip 6pt

For $x\in S$ let $\Ux\subset\Xx\supset\Zx$ denote the fibers of
$\cU\subset\cX\supset\cZ$.

\Def 2-1. \it Let $G_j\in H^0(\cX,\cO_{\cX}(e_j))$ be a non-zero
element defining $\cZ_i\subset \cX$.
\roster
\item
For $1\leq j\leq s-1$ put
$$g_j=(G_j^{e_s}/G_s^{e_j})\in \Gamma(\cU,\cO_{\cU}^*)=CH^1(\cU,1)$$
where $CH^i(*,j)$ is the Bloch's higher Chow group (cf. [Bl]).
\item
Define the index set
$J=\{\sigma=(j_1,\dots,j_m)|\spa 1\leq j_1<\cdots<j_m\leq s-1\}$
where $m=n-r$ is the relative dimension of $\cX/S$.
For $\sigma=(j_1,\dots,j_m)\in J$ let
$$ \ccU=\{g_{j_1},\dots,g_{j_m}\}\in CH^m(\cU,m)$$
be defined by the product
$CH^1(\cU,1)\otimes\cdots\otimes CH^1(\cU,1)\to CH^m(\cU,m)$.
For $x\in S$ let $\cUx\in CH^m(\Ux,m)$ be the restriction of $\ccU$.
\item
For $\sigma=(j_1,\dots,j_m)\in J$ let
$$ \wcU=\dlog{g_{j_1}}\wedge \cdots \wedge \dlog{g_{j_m}}\in
\Gamma(S,f_* \WcXZ m)$$
where $\WcXZ p=\os{p}{\wedge}\WcXZ 1$ with $\WcXZ 1$, the sheaf of
relative differentials on $\cX$ over $S$ with logarithmic poles along
$\cZ$.
\endroster
\vskip 6pt\rm

\Lem 2-1. \it
$\wcU$ with $\sigma\in J$ are linearly independent over $\cO_S$.
\rm\demo{Proof}
For $\tau=(j_1,\dots,j_m)\in J$ write
$Z_{\tau}=Z_{j_1}\cap\cdots\cap Z_{j_m}$. By taking Poincar\'e
residues along $Z_{j_1},\dots Z_{j_m}$, we get the $\cO_S$-linear map
$$ Res_{\tau}\scs \Gamma(S,f_* \WcXZ m)\to
\Gamma(Z_{\tau},\cO_{Z_{\tau}})$$
and we have
$$ Res_{\tau}(\wcU)=\left\{\aligned
& 1\quad\text{ if } \sigma=\tau,\\
& 0\quad\text{ if } \sigma\not=\tau.\\
\endaligned\right.$$
This proves Lem.(2-1).
\qed\enddemo
\vskip 5pt


\Def 2-2. \it For $x\in S$ let
$$ \roxlog \scs \TxS \to H^1(\Xx,\TXZx)
\quad (resp. \spa \cXZx: B_1(0) \to H^1(\Xx,\TXZx))$$
be the Kodaira-Spencer map (resp. the map in Th.(I)(2) for $(\Xx,\Zx)$).
We define
$$\ccXZS=\underset{\spa x\in S}\to{\max}
\{\dim_{k}(\Im(\cXZx)/\Im(\cXZx)\cap \Im(\roxlog))\}.$$
\vskip 5pt\rm

\Rem 2-1. If $n-r\geq 2$ and $\Xx$ is not a $K3$ surface,
$\cXZx$ is surjective so that
$$\ccXZS=\underset{x\in S}\to{\max}\{\dim_{k(x)}
(\Coker(\rocXZ)\otimes_{\cO_S} k(x))\} .$$

\Rem 2-2. We may use the following more intuitive invariant than $\ccXZS$.
Let $P^d\subset P$ be as in \S1. The dual projective space
$$\Pd(P^d)=\Bbb P_k^{N(n,d)}\quad (N(n,d)=\binom{n+d}{d}-1)$$
parametrizes hypersurfaces
$Y\subset \Bbb P^n$ of degree $d$ defined over $k$. Let
$$ B\subset \underset{1\leq i\leq r}\to{\prod} \Bbb P_k^{N(n,d_i)} \times
\underset{1\leq j\leq s}\to{\prod} \Bbb P_k^{N(n,e_j)}$$
be the Zariski open subset parametrizing such $(Y_\nu)_{1\leq \nu\leq r+s}$
that $Y_1+\cdots+Y_{r+s}$ is a simple normal crossing divisor on $\Bbb P^n_k$.
We consider the family
$$\cXB \hookleftarrow \cZB=\underset{1\leq j\leq s}\to{\cup}
\cZBj \quad\text{ over } B
$$
whose fibers are
$ X \hookleftarrow Z=\underset{1\leq j\leq s}\to{\cup}Z_j $
with $X=Y_1\cap\cdots \cap Y_r$ and $Z_j=X\cap Y_{r+j}$.
Let $T\subset B$ be a non-singular locally closed subvariety of codimension
$c\geq 0$ and let $S\to T$ be a dominant map. Assume that the family (2-1)
is the pullback of $(\cXB,\cZB)/B$ via $S\to B$.
Then we have $\ccXZS\leq c$ and the statement of Th.(2-1) holds with
$\ccXZS$ replaced by $c$.
\vskip 6pt

Now we prove Th.(2-1).
We only show the second assertion and leave the first to the readers.
The fact $\wcU\in \Gamma(S,\Ker(\nabmo))$ follows from the fact that
$\wcU$ lies in the image of
$$H^0(\cX,\WcXkZ m)\to \Gamma(S,f_*\WcXZ m),$$
where $\WcXkZ \cdot$ is the sheaf of differential forms of $\cX$ over $k$
with logarithmic poles along $\cZ$.
Fix $0\in S$ and let $X\supset Z$
be the fibers of the family (2-1). Let $\Sigma\subset H^0(X,\WXZ m)$
be the subspace generated by $\wcU(0)$ with $\sigma\in J$.
It suffices to show the injectivity of
$$H^0(X,\WXZ m)/\Sigma \to \Omega_{S,0}^1\otimes H^1(X,\WXZ {m-1})$$
that is induced by $\nabmo$.
By Th.(II) and Th.(II') in \S1 this is reduced to show the surjectivity of
$$ V\otimes B_{n-r-1}(\d-n-1) \to B_{n-r}(\d-n-1)$$
where $B_1(0)\supset V:=\cXZ^{-1}(\Im(T_0(S) @>{\roolog}>> H^1(X,\TXZ)))$
(cf. Def.(2-2)). By definition $V$ is of codimension$\leq c$ in $B_1(0)$.
Hence the desired assertion follows from Th.(III)$(i)$ in \S1.
\qed
\vskip 5pt

\vskip 20pt

\input amstex
\documentstyle{amsppt}
\hsize=16cm
\vsize=23cm

\def\Th#1.{\vskip 6pt \medbreak\noindent{\bf Theorem(#1).}}
\def\Cor#1.{\vskip 6pt \medbreak\noindent{\bf Cororally(#1).}}
\def\Conj#1.{\vskip 6pt \medbreak\noindent{\bf Conjecture(#1).}}
\def\Pr#1.{\vskip 6pt \medbreak\noindent{\bf Proposition(#1).}}
\def\Lem#1.{\vskip 6pt \medbreak\noindent{\bf Lemma(#1).}}
\def\Rem#1.{\vskip 6pt \medbreak\noindent{\it Remark(#1).}}
\def\Fact#1.{\vskip 6pt \medbreak\noindent{\it Fact(#1).}}
\def\Claim#1.{\vskip 6pt \medbreak\noindent{\it Claim(#1).}}
\def\Def#1.{\vskip 6pt \medbreak\noindent{\bf Definition\bf(#1)\rm.}}

\def\qwith{\quad\hbox{with }}
\def\mathrm#1{\rm#1}

\def\isom{@>\cong>>}
\def\Spec{{\operatorname{Spec}}}

\def\dim{{\operatorname{dim}}}

\def\Coker{{\text{\rm Coker}}}
\def\dim{\hbox{\rm dim}}

\def\Im{\hbox{\rm Im}}
\def\Ker{\hbox{\rm Ker}}
\def\Coker{\hbox{\rm Coker}}
\def\min{\hbox{\rm min}}

\def\Gal{\hbox{\mathrm{Gal}}}
\def\GL{\hbox{\mathrm{GL}}}

\def\P{{\Bbb{P}}}
\def\bP{{\Bbb{P}}}

\def\cHom{{\Cal{H}}om}

\def\cF{{\Cal{F}}}
\def\cO{{\Cal{O}}}
\def\cX{{\Cal{X}}}
\def\cXB{{\Cal{X}_B}}

\def\cZ{{\Cal{Z}}}
\def\cZB{{\Cal{Z}_B}}
\def\cZBj{{\Cal{Z}_{B,j}}}

\def\wtd#1{\widetilde{#1}}

\def\l{\ell}
\def\d{{\bold{d}}}
\def\e{{\bold{e}}}

\def\ra{\rightarrow}

\def\ot{\otimes}


\def\us#1#2{\underset{#1}\to{#2}}
\def\os#1#2{\overset{#1}\to{#2}}

\def\qaq{\quad\hbox{ and }\quad}
\def\qfor{\quad\hbox{ for }}

\def\spa{\hbox{ }}
\def\scs{\spa : \spa}

\def\sp{\hbox{}}

\def\ul#1{\underline{#1}}

\def\ba{\bold{a}}
\def\bb{\bold{b}}
\def\ua{\underline{a}}
\def\ub{\underline{b}}

\def\ue{\underline{e}}
\def\ud{\underline{d}}

\def\md{\delta_{\text{min}}}

\def\WS#1{\Omega_S^{#1}}

\def\otO{\otimes_{\cO}}
\def\cU{{\Cal U}}

\def\onab{\overline{\nabla}}

\def\HO{H_{\cO}}
\def\HQ{H_{\Bbb Q}}
\def\HC{H_{\Bbb C}}
\def\HQcU#1#2{\HQ^{#1}(\cU/S)(#2)}
\def\HCcU#1{\HC^{#1}(\cU/S)}
\def\HOcU#1{\HO^{#1}(\cU/S)}

\def\HcU#1#2{H^{#1,#2}(\cU/S)}

\def\HcU#1#2{H^{#1,#2}(\cU/S)}
\def\HcXp#1#2{H^{#1,#2}(\cX/S)_{prim}}

\def\HcUan#1#2{H^{#1,#2}(\cU/S)^{an}}

\def\HQcU#1#2{\HQ^{#1}(\cU/S)(#2)}

\def\HQcU#1#2{\HQ^{#1}(\cU/S)(#2)}

\def\HCcU#1{\HC^{#1}(\cU/S)}
\def\HCcX#1{\HC^{#1}(\cX/S)}
\def\HCcU#1{\HC^{#1}(\cU/S)}
\def\HCcXp#1{\HC^{#1}(\cX/S)_{prim}}

\def\HOcU#1{\HO^{#1}(\cU/S)}
\def\HOcX#1{\HO^{#1}(\cX/S)}
\def\HOcU#1{\HO^{#1}(\cU/S)}

\def\WSan#1{\Omega_{\San}^{#1}}

\def\Xx{X_x}
\def\Zx{Z_x}
\def\Ux{U_x}
\def\xb{\overline{x}}
\def\Uxb{U_{\xb}}

\def\Zxx{Z_{x}}
\def\Zst{Z}
\def\cZst{\cZ}

\def\Yst{Y}

\def\rocXZ{\kappa_{(\cX,\cZ)}}

\def\roolog{\kappa_0^{log}}

\def\roxlog{\kappa_x^{log}}

\def\TS{\Theta_S}
\def\TxS{T_x S}

\def\UC{U_{\Bbb C}}

\def\SC{S_{\Bbb C}}

\def\WS#1{\Omega_S^{#1}}
\def\WX#1{\Omega_X^{#1}}

\def\TX{T_X}

\def\TXZx{T_{\Xx}(-\log \Zxx)}

\def\WXZ#1{\Omega_X^{#1}(\log \Zst)}

\def\WcXZ#1{\Omega_{\cX/S}^{#1}(\log \cZst)}

\def\WcXkZ#1{\Omega_{\cX/k}^{#1}(\log \cZst)}

\def\WPnY#1{\Omega_{\Bbb P^n}^{#1}(\log \Yst)}
\def\TXZ{T_X(-\log \Zst)}
\def\TcXZS{T_{\cX/S}(-\log \cZst)}

\def\fXZ#1{\phi_{X,Z}^{#1}}
\def\cXZ{\psi_{(X,Z)}}

\def\cXZx{\psi_{(\Xx,\Zx)}}

\def\eXZ{\eta_{(X,Z)}}

\def\ccXZS{c_S(\cX,\cZ)}
\def\ccXZtS{c_{\wtd{S}}(\wtd{\cX},\wtd{\cZ})}

\def\tS{\wtd{S}}

\def\scs{\hbox{ }:\hbox{ }}
\def\onab{\overline{\nabla}}

\def\dlog#1{\frac{d#1}{#1}}
\def\chB#1#2#3{ch_{B,#3}^{#1,#2}}
\def\chD#1#2#3{ch_{D,#3}^{#1,#2}}
\def\chet#1#2#3{ch_{et,#3}^{#1,#2}}
\def\chcont#1#2#3{ch_{cont,#3}^{#1,#2}}
\def\Ql{\Bbb Q_{\ell}}
\def\Zl{\Bbb Z_{\ell}}
\def\etab{\overline{\eta}}

\def\Pd{\overset{\vee}\to{\Bbb P}}
\def\kb{\overline{k}}
\def\xb{\overline{x}}

\def\pitop{\pi_1^{top}}
\def\pialg{\pi_1^{alg}}

\def\wcU{\omega_{\cU/S}(\sigma)}
\def\dcU{\delta_{\cU/S}(\sigma)}

\def\ccU{c_{\cU/S}(\sigma)}
\def\cUx{c_{\Ux}(\sigma)}
\def\cUo{c_{U}(\sigma)}
\def\cUetab{c_{U_{\etab}}(\sigma)}

\def\Dqi{\Delta_{\gamma_i}}
\def\Dq{\Delta_{\gamma}}
\def\spa{\hbox{ }}

\def\nabmo{\overline{\nabla}^{m,0}}
\def\nab#1#2{\overline{\nabla}^{#1,#2}}

\def\Pol{P}
\def\h#1#2{h_{#1}(#2)}

\def\San{S_{an}}
\def\Szar{S_{zar}}

\head \S3. Beilinson's Hodge conjecture. \endhead
\vskip 8pt

In this section we assume that $k=\Bbb C$.
Let $\San$ be the analytic site on $S(\Bbb C)$.
For a coherent sheaf $\cF$ on $\Szar$ let $\cF^{an}$ be the associated analytic
sheaf on $\San$. We introduce local systems on $\San$
$$\HQcU q p=R^q g_*\Bbb Q(p) \qaq \HCcU q =R^q g_*\Bbb C,$$
where $g:\cU\to S$ is the natural morphism.
Let $\HOcU q$ be the sheaf of holomorphic sections of $\HCcU q$
and let $F^p \HOcU q\subset \HOcU q$ be the holomorphic subbundle given
by the Hodge filtration on the cohomology of fibers of $\cU/S$.
We have the analytic Gauss-Manin connection
$$ \nabla\scs \HOcU q \to \WSan 1 \otimes \HOcU q \qwith
\Ker(\nabla)=\HCcU q$$
that satisfies
$\nabla(F^p\HOcU q) \subset \WSan 1 \otimes F^{p-1}\HOcU q.$
The induced map
$$ F^p\HOcU {p+q}/F^{p+1} \to \WSan 1 \otimes F^{p-1}\HOcU {p+q}/F^p$$
is identified with $(\nab p q)^{an}$ via the identification
$F^p \HOcU {p+q}/F^{p+1} = \HcUan p q$.
Therefore Th.(2-1) implies the following.

\Th 3-1. \it Assume $\d\geq n+1+\ccXZS$.
\roster
\item
If $s\geq 1$, $F^1\HOcU m \cap \HCcU m$
is generated over $\Bbb C$ by
$$\dcU:=[\frac{d g_{j_1}}{g_{j_1}}]\cup \cdots \cup [\frac{d g_{j_m}}{g_{j_m}}]
\in \Gamma(\San,\HQcU m m) \qwith \sigma=(j_1,\dots,j_m)\in J,$$
where
$[\frac{d g_{j}}{g_{j}}]\in \Gamma(\San,\HQcU 1 1)$ is the cohomology class of
$\frac{d g_{j}}{g_{j}}$.
\item
If $s=0$, $F^1\HOcX m \cap \HCcXp m=0$, where $\HCcXp q$ is the primitive part
of $\HCcX q$.
\endroster
\rm
\vskip 5pt

\Th 3-2. \it Assume
$\md(n-r-1)+\d\geq n+1+\ccXZS$.
There exists $E\subset S(\Bbb C)$ such that:
\roster
\item"$(i)$"
$E$ is the union of countable many closed analytic subset of
codimension$\geq 1$.
\item"$(ii)$"
For $\forall x\in S(\Bbb C)-E$,
$$H_B^m(\Ux(\Bbb C),\Bbb Q(m))\cap F^m H_B^m(\Ux(\Bbb C),\Bbb C) =
\underset{\sigma\in J}\to{\bigoplus}\spa \Bbb Q\cdot \chB m m {\Ux}(\cUx).$$
\endroster
\rm\demo{Proof}
Fix a base point $0\in S$ and let $U\subset X\supset Z$ be the fibers of
$\cU\subset \cX \supset \cZ$ over $0$.
Let $0\in \Delta\subset S(\Bbb C)$ be an open disk.
Identifying $H_B^m(U,\Bbb Q(m))$ with $\Gamma(\Delta,\HQcU m m)$, we set
$$ \Dq=\{x\in \Delta|\spa \gamma(x)\in F^m H_B^m(\Ux,\Bbb C)\}
\qfor \gamma\in H_B^m(U,\Bbb Q(m)).$$
Let $\wtd {\gamma}\in \Gamma(\Delta,\HOcU m/F^m\HOcU m)$ be the image
of $\gamma\in \Gamma(\Delta,\HQcU m m)\subset \Gamma(\Delta,\HOcU m).$
Then $\Dq\subset \Delta$ is the zero locus of $\wtd \gamma$ and hence it is
an analytic subset. Write
$H_B^m(U,\Bbb Q(m))=\{\gamma_i\}_{i\in I}$ as a set.
Note that $I$ is a countable set. Setting
$A=\{i\in I|\spa \Dqi=\Delta\}$ and
$B=\{i\in I|\spa \Dqi\subsetneq \Delta\}$, we have
$I=A\cup B$ and $A\cap B=\emptyset$.
We put $E_\Delta=\underset{i\in B}\to{\cup} \Dqi$.
It suffices to show Th.(3-2)$(ii)$ holds for
$\forall x\in \Delta-E_\Delta$.
By definition, for $\forall x\in \Delta-E_\Delta$ we have as a set
$$H_B^m(\Ux,\Bbb Q(m))\cap F^m H_B^m(\Ux,\Bbb C)=\{\gamma_i(x)\}_{i\in A}$$
that implies
$$ \aligned
 H_B^m(\Ux,\Bbb Q(m))\cap F^m H_B^m(\Ux,\Bbb C)
&\simeq \Gamma(\Delta,\HQcU m m\cap F^m\HOcU m)\\
&\subset \Gamma(\Delta,\Ker(F^m\HOcU m @>\nabla>>\WS 1\otO F^{m-1}\HOcU m)\\
\endaligned
\leqno(*)$$
Note that $\chB m m {\Ux}(\cUx)=\wcU(x)$ under the natural identification of
$ F^m H_B^m(\Ux,\Bbb C)$ with the fiber of $\HcU m 0$ over $x$.
Hence the desired assertion follows from Th.(2-1) and $(*)$ and Lem.(2-1).
\qed
\enddemo
\vskip 6pt

\vskip 20pt

\input amstex
\documentstyle{amsppt}
\hsize=16cm
\vsize=23cm

\def\Th#1.{\vskip 6pt \medbreak\noindent{\bf Theorem(#1).}}
\def\Cor#1.{\vskip 6pt \medbreak\noindent{\bf Cororally(#1).}}
\def\Conj#1.{\vskip 6pt \medbreak\noindent{\bf Conjecture(#1).}}
\def\Pr#1.{\vskip 6pt \medbreak\noindent{\bf Proposition(#1).}}
\def\Lem#1.{\vskip 6pt \medbreak\noindent{\bf Lemma(#1).}}
\def\Rem#1.{\vskip 6pt \medbreak\noindent{\it Remark(#1).}}
\def\Fact#1.{\vskip 6pt \medbreak\noindent{\it Fact(#1).}}
\def\Claim#1.{\vskip 6pt \medbreak\noindent{\it Claim(#1).}}
\def\Def#1.{\vskip 6pt \medbreak\noindent{\bf Definition\bf(#1)\rm.}}

\def\qwith{\quad\hbox{with }}
\def\mathrm#1{\rm#1}

\def\isom{@>\cong>>}
\def\Spec{{\operatorname{Spec}}}

\def\dim{{\operatorname{dim}}}

\def\Coker{{\text{\rm Coker}}}
\def\dim{\hbox{\rm dim}}

\def\Im{\hbox{\rm Im}}
\def\Ker{\hbox{\rm Ker}}
\def\Coker{\hbox{\rm Coker}}
\def\min{\hbox{\rm min}}

\def\Gal{\hbox{\mathrm{Gal}}}
\def\GL{\hbox{\mathrm{GL}}}

\def\P{{\Bbb{P}}}
\def\bP{{\Bbb{P}}}

\def\cHom{{\Cal{H}}om}

\def\cF{{\Cal{F}}}
\def\cO{{\Cal{O}}}
\def\cX{{\Cal{X}}}
\def\cXB{{\Cal{X}_B}}

\def\cZ{{\Cal{Z}}}
\def\cZB{{\Cal{Z}_B}}
\def\cZBj{{\Cal{Z}_{B,j}}}

\def\wtd#1{\widetilde{#1}}

\def\l{\ell}
\def\d{{\bold{d}}}
\def\e{{\bold{e}}}

\def\ra{\rightarrow}

\def\ot{\otimes}


\def\us#1#2{\underset{#1}\to{#2}}
\def\os#1#2{\overset{#1}\to{#2}}

\def\qaq{\quad\hbox{ and }\quad}
\def\qfor{\quad\hbox{ for }}

\def\spa{\hbox{ }}
\def\scs{\spa : \spa}

\def\sp{\hbox{}}

\def\ul#1{\underline{#1}}

\def\ba{\bold{a}}
\def\bb{\bold{b}}
\def\ua{\underline{a}}
\def\ub{\underline{b}}

\def\ue{\underline{e}}
\def\ud{\underline{d}}

\def\md{\delta_{\text{min}}}

\def\WS#1{\Omega_S^{#1}}

\def\otO{\otimes_{\cO}}
\def\cU{{\Cal U}}

\def\onab{\overline{\nabla}}

\def\HO{H_{\cO}}
\def\HQ{H_{\Bbb Q}}
\def\HC{H_{\Bbb C}}
\def\HQcU#1#2{\HQ^{#1}(\cU/S)(#2)}
\def\HCcU#1{\HC^{#1}(\cU/S)}
\def\HOcU#1{\HO^{#1}(\cU/S)}

\def\HcU#1#2{H^{#1,#2}(\cU/S)}

\def\HcU#1#2{H^{#1,#2}(\cU/S)}
\def\HcXp#1#2{H^{#1,#2}(\cX/S)_{prim}}

\def\HcUan#1#2{H^{#1,#2}(\cU/S)^{an}}

\def\HQcU#1#2{\HQ^{#1}(\cU/S)(#2)}

\def\HQcU#1#2{\HQ^{#1}(\cU/S)(#2)}

\def\HCcU#1{\HC^{#1}(\cU/S)}
\def\HCcX#1{\HC^{#1}(\cX/S)}
\def\HCcU#1{\HC^{#1}(\cU/S)}
\def\HCcXp#1{\HC^{#1}(\cX/S)_{prim}}

\def\HOcU#1{\HO^{#1}(\cU/S)}
\def\HOcX#1{\HO^{#1}(\cX/S)}
\def\HOcU#1{\HO^{#1}(\cU/S)}

\def\WSan#1{\Omega_{\San}^{#1}}

\def\Xx{X_x}
\def\Zx{Z_x}
\def\Ux{U_x}
\def\xb{\overline{x}}
\def\Uxb{U_{\xb}}

\def\Zxx{Z_{x}}
\def\Zst{Z}
\def\cZst{\cZ}

\def\Yst{Y}

\def\rocXZ{\kappa_{(\cX,\cZ)}}

\def\roolog{\kappa_0^{log}}

\def\roxlog{\kappa_x^{log}}

\def\TS{\Theta_S}
\def\TxS{T_x S}

\def\UC{U_{\Bbb C}}

\def\SC{S_{\Bbb C}}

\def\WS#1{\Omega_S^{#1}}
\def\WX#1{\Omega_X^{#1}}

\def\TX{T_X}

\def\TXZx{T_{\Xx}(-\log \Zxx)}

\def\WXZ#1{\Omega_X^{#1}(\log \Zst)}

\def\WcXZ#1{\Omega_{\cX/S}^{#1}(\log \cZst)}

\def\WcXkZ#1{\Omega_{\cX/k}^{#1}(\log \cZst)}

\def\WPnY#1{\Omega_{\Bbb P^n}^{#1}(\log \Yst)}
\def\TXZ{T_X(-\log \Zst)}
\def\TcXZS{T_{\cX/S}(-\log \cZst)}

\def\fXZ#1{\phi_{X,Z}^{#1}}
\def\cXZ{\psi_{(X,Z)}}

\def\cXZx{\psi_{(\Xx,\Zx)}}

\def\eXZ{\eta_{(X,Z)}}

\def\ccXZS{c_S(\cX,\cZ)}
\def\ccXZtS{c_{\wtd{S}}(\wtd{\cX},\wtd{\cZ})}

\def\tS{\wtd{S}}

\def\scs{\hbox{ }:\hbox{ }}
\def\onab{\overline{\nabla}}

\def\dlog#1{\frac{d#1}{#1}}
\def\chB#1#2#3{ch_{B,#3}^{#1,#2}}
\def\chD#1#2#3{ch_{D,#3}^{#1,#2}}
\def\chet#1#2#3{ch_{et,#3}^{#1,#2}}
\def\chcont#1#2#3{ch_{cont,#3}^{#1,#2}}
\def\Ql{\Bbb Q_{\ell}}
\def\Zl{\Bbb Z_{\ell}}
\def\etab{\overline{\eta}}

\def\Pd{\overset{\vee}\to{\Bbb P}}
\def\kb{\overline{k}}
\def\xb{\overline{x}}

\def\pitop{\pi_1^{top}}
\def\pialg{\pi_1^{alg}}

\def\wcU{\omega_{\cU/S}(\sigma)}
\def\dcU{\delta_{\cU/S}(\sigma)}

\def\ccU{c_{\cU/S}(\sigma)}
\def\cUx{c_{\Ux}(\sigma)}
\def\cUo{c_{U}(\sigma)}
\def\cUetab{c_{U_{\etab}}(\sigma)}

\def\Dqi{\Delta_{\gamma_i}}
\def\Dq{\Delta_{\gamma}}
\def\spa{\hbox{ }}

\def\nabmo{\overline{\nabla}^{m,0}}
\def\nab#1#2{\overline{\nabla}^{#1,#2}}

\def\Pol{P}
\def\h#1#2{h_{#1}(#2)}

\def\San{S_{an}}
\def\Szar{S_{zar}}

\head \S4. Beilinson's Tate conjecture. \endhead
\vskip 8pt

In this section we show Th.(0-1)(2) in the introduction.
It will follows from Th.(4-3) below by using the theory of Hilbert set
(cf. [La]). Let the assumption be as in the beginning of \S2.
Write $m=n-r$.
\vskip 6pt

\Th 4-1. \it Assume $\d \geq n+1+\ccXZS$.
Let $\etab$ be a geometric generic point of $S$ and
let $U_{\etab}\subset X_{\etab}$ be the fibers of $\cU\subset\cX$ over $\etab$.
\roster
\item
Assuming $s\geq 1$, we have
$$H^m_{et}(U_{\etab},\Ql(m))^{\pi_1(S,\etab)}=
\underset{\sigma\in J}\to{\bigoplus}\spa
\Ql\cdot \chet m m {U_{\etab}}(\cUetab).$$
\item
$H^m_{et}(X_{\etab},\Ql(m))_{prim}^{\pi_1(S,\etab)}=0.$
\endroster
\rm\vskip 5pt

\Rem 4-1. \it Let $\wtd S \to S$ be etale and
$\wtd{\cX}\supset \wtd{\cZ}$ be the base change. Then $\ccXZS=\ccXZtS$
Hence Th.(4-1) holds if one replaces $\pi_1(S,\etab)$ by any open subgroup of
finite index.
\rm\vskip 5pt

\Th 4-2. \it Assume $\d \geq n+1+\ccXZS$.
Let $U\subset X$ be the fibers of $\cU\subset\cX$ over a fixed base point
$0\in S(\Bbb C)$.
\roster
\item
Assuming $s\geq 1$, we have
$$H_B^m(U(\Bbb C),\Bbb Q(m))^{\pi_1(S,0)}=
\underset{\sigma\in J}\to{\bigoplus}\spa \Bbb Q\cdot\chB m m U(\cUo).$$
\item
$H_B^m(X(\Bbb C),\Bbb Q(m))_{prim}^{\pi_1(S,0)}=0.$
\endroster \rm
\vskip 4pt

First we deduce Th.(4-1) from Th.(4-2).
By the Lefschetz principle we may assume that $k$ is a subfield of $\Bbb C$
finitely generated over $\Bbb Q$. We fix an embedding
$k(\etab)\hookrightarrow \Bbb C$ and let $0\in S(\Bbb C)$ be the corresponding
$\Bbb C$-valued point of $S$. Write $\SC=S\otimes_k \Bbb C$ and put
$\UC=\cX\times_S \Spec(\Bbb C)$ via $\Spec(\Bbb C)@>0>> S$.
We have the comparison isomorphisms ([SGA4, XVI Th.4.1])
$$ H_B^m(\UC,\Bbb Q)\otimes \Ql \isom H^m_{et}(\UC,\Ql)\isom
 H^m_{et}(U_{\etab},\Ql)$$
that are equivariant with respect to the maps of topological and algebraic
fundamental groups:
$$ \pitop(\SC,0) \to \pialg(\SC,0) \to \pialg(S,\etab).$$
The desired assertion follows at once from this.
\qed

\vskip 6pt

First we show Th.(4-2)(2). Recall the notation in Th.(3-1).
Write $H=H_B^{m}(X(\Bbb C),\Bbb Q)_{prim}^{\pi_1(S,0)}$. We have the natural
isomorphism
$$ H\otimes\Bbb C\isom \Gamma(\San,\HCcXp {m}).$$
By [D] $H$ is a sub-Hodge structure of $H_B^{m}(X(\Bbb C),\Bbb Q)$ so that
we have the Hodge decomposition
$$H\otimes\Bbb C=\underset{p+q=m}\to{\bigoplus} H^{p,q} \qwith
H^{p,q}\subset \Gamma(\San,F^p\HOcX {m} \cap \HCcXp {m}).$$
By Th.(3-1) this implies $H^{p,q}=0$ for $p\geq 1$ which implies $H=0$
by the Hodge symmetry.
\qed
\vskip 4pt

In order to show Th.(4-2)(1) we need the following

\Lem 4-1. \it Let the assumption be as in Th.(4-2). Then
$H:=H_B^m(U(\Bbb C),\Bbb Q)^{\pi_1(S,0)}$ is a submixed Hodge structure of
$H_B^m(U(\Bbb C),\Bbb Q)$ and $F^1\HC=\HC$ where $\HC=H\otimes\Bbb C$.
\rm\demo{Proof}
The fact that $H$ is a submixed Hodge structure follows from the theory
of mixed Hodge modules [SaM]. We show the second assertion.
We recall that the graded subquotients of the weight filtration of
$H_B^m(U(\Bbb C),\Bbb Q)$ is given by
$$ Gr^W_{m+p}H^m(U,\Bbb Q)=
\underset{1\leq j_1<\cdots<j_p\leq s}\to{\bigoplus}
H^{m-p}(Z_{j_1}\cap\cdots\cap Z_{j_p},\Bbb Q(-p))_{prim}.
\quad (0\leq p\leq m)$$
Therefore it suffices to note that
$ H_B^{m}(X(\Bbb C),\Bbb Q)_{prim}^{\pi_1(S,0)}=0$ as we have just shown.
\qed
\enddemo
\vskip 4pt

Finally we complete the proof of Th.(4-2)(1).
The linear independence of $\chB m m U(\cUo)$ with $\sigma\in J$ follows from
Lem.(2-1). We show that they span $H_B^m(U(\Bbb C),\Bbb Q(m))^{\pi_1(S,0)}$.
Recall the notation in Th.(3-1). By Lem.(4-1) we have
$$ H_B^m(U(\Bbb C),\Bbb Q(m))^{\pi_1(S,0)} \simeq \Gamma(\San,\HQcU m m)
\subset \Gamma(\San,\HCcU m \cap F^1\HOcU m).$$
Hence the assertion follows from Th.(3-1)(1) by noting
$\chB m m {U}(\cUo)=\dcU(0)$.
\qed
\vskip 6pt

\def\tpi{\widetilde{\pi}}

\Th 4-3. \it Let the assumption be as in Th.(0-1)(2).
Then there exists an irreducible variety $\tS$ over $k$
with a finite etale covering $\phi:\tS \to S$ for which the following holds:
Let $\tpi=\pi\circ\phi$ and let $H\subset \Bbb P^N(k)$ be the subset of such
points that $\tpi^{-1}(y)$ is irreducible. For $\forall x\in S$ such that
$\pi(x)\in H$ we have
$$H^m_{et}(\Uxb,\Ql(m))^{Gal(k(\xb)/k(x))}=
\underset{\sigma\in J}\to{\bigoplus}\spa\Ql \cdot\chet m m {\Ux}(\cUx).$$
\rm\demo{Proof} (cf. [T] and [BE])
By choosing a $k$-rational point $0$ of $S$, we get the decomposition
$$\pi_1(S,\etab)=\pi_1(S,\etab)^{geo}\times \Gal(\kb/k),$$
where $\pi_1(S,\etab)^{geo}$ classifies the finite etale coverings of $S$
that completely decompose over $0$. Let
$$\Gamma=\Im\big(\pi_1(S,\etab)^{geo} \to
\GL_{\Ql}(H^m_{et}(U_{\etab},\Ql(m))\big).$$
From the fact that $\Gamma$ contains an $\ell$-adic Lie group as a subgroup
of finite index, we have the following fact (cf. [T]):
There exists a subgroup $\Gamma'\subset \Gamma$ of finite index such that
a continuous homomorphism of pro-finite group $G \to \Gamma$ is surjective
if and only if $G\to\Gamma\to \Gamma/\Gamma'$ is surjective as a map of sets.
Let $\phi:\tS\to S$ be a finite etale covering that corresponds to the inverse
image of $\Gamma'$ in $\pi_1(S,\etab)^{geo}$ and let $H$ be defined as in
Th.(4-3). Fix $x\in S$ with $\pi(x)\in H$ and let $\xb$ be a geometric
point of $x$. By choosing a ``path" $\xb\to \etab$, we get the map
$\Gal(k(\xb)/k(x))=\pi_1(x,\xb)@>\iota>> \pi_1(S,\etab)$.
By the definition of $H$ the composite of $\iota$ with
$\pi_1(S,\etab)\to\Gamma\to\Gamma/\Gamma'$ is surjective so that
$\Gal(k(\xb)/k(x))$ surjects onto $\Gamma$ by the above fact.
This implies the isomorphism
$$ H^m_{et}(U_{\xb},\Ql(m))^{\pi_1(x,\xb)}\isom
H^m_{et}(U_{\etab},\Ql(m))^{\pi_1(S,\etab)}.$$
Now Th.(4-3) follows from Th.(4-1).
\qed
\enddemo
\vskip 5pt

Now Th.(0-1)(2) is a consequence of Th.(4-3) and the following result
(cf. [La]):

\Lem 4-2. \it Let $k$ be a number field and let $V$ be an irreducible
variety over $k$. Let $\pi:V\to \Bbb P_k^N$ be an etale morphism
and let $H\subset \Bbb P_k^N(k)$ be the subset of such points $x$
that $\pi^{-1}(x)$ is irreducible.
Let $\Sigma$ be any finite set of primes of $k$ and let $k_v$ be the
completion of $k$ at $v$. Then the image of $H$ in
$\prod_{v\in \Sigma} \Bbb P_k^N(k_v)$ is dense.
\rm

\vskip 20pt

\input amstex
\documentstyle{amsppt}
\hsize=16cm
\vsize=23cm

\def\Th#1.{\vskip 6pt \medbreak\noindent{\bf Theorem(#1).}}
\def\Cor#1.{\vskip 6pt \medbreak\noindent{\bf Cororally(#1).}}
\def\Conj#1.{\vskip 6pt \medbreak\noindent{\bf Conjecture(#1).}}
\def\Pr#1.{\vskip 6pt \medbreak\noindent{\bf Proposition(#1).}}
\def\Lem#1.{\vskip 6pt \medbreak\noindent{\bf Lemma(#1).}}
\def\Rem#1.{\vskip 6pt \medbreak\noindent{\it Remark(#1).}}
\def\Fact#1.{\vskip 6pt \medbreak\noindent{\it Fact(#1).}}
\def\Claim#1.{\vskip 6pt \medbreak\noindent{\it Claim(#1).}}
\def\Def#1.{\vskip 6pt \medbreak\noindent{\bf Definition\bf(#1)\rm.}}

\def\qwith{\quad\hbox{with }}
\def\mathrm#1{\rm#1}

\def\isom{@>\cong>>}
\def\Spec{{\operatorname{Spec}}}

\def\dim{{\operatorname{dim}}}

\def\Coker{{\text{\rm Coker}}}
\def\dim{\hbox{\rm dim}}

\def\Im{\hbox{\rm Im}}
\def\Ker{\hbox{\rm Ker}}
\def\Coker{\hbox{\rm Coker}}
\def\min{\hbox{\rm min}}

\def\Gal{\hbox{\mathrm{Gal}}}
\def\GL{\hbox{\mathrm{GL}}}

\def\P{{\Bbb{P}}}
\def\bP{{\Bbb{P}}}

\def\cHom{{\Cal{H}}om}

\def\cF{{\Cal{F}}}
\def\cO{{\Cal{O}}}
\def\cX{{\Cal{X}}}
\def\cXB{{\Cal{X}_B}}

\def\cZ{{\Cal{Z}}}
\def\cZB{{\Cal{Z}_B}}
\def\cZBj{{\Cal{Z}_{B,j}}}

\def\wtd#1{\widetilde{#1}}

\def\l{\ell}
\def\d{{\bold{d}}}
\def\e{{\bold{e}}}

\def\ra{\rightarrow}

\def\ot{\otimes}


\def\us#1#2{\underset{#1}\to{#2}}
\def\os#1#2{\overset{#1}\to{#2}}

\def\qaq{\quad\hbox{ and }\quad}
\def\qfor{\quad\hbox{ for }}

\def\spa{\hbox{ }}
\def\scs{\spa : \spa}

\def\sp{\hbox{}}

\def\ul#1{\underline{#1}}

\def\ba{\bold{a}}
\def\bb{\bold{b}}
\def\ua{\underline{a}}
\def\ub{\underline{b}}

\def\ue{\underline{e}}
\def\ud{\underline{d}}

\def\md{\delta_{\text{min}}}

\def\WS#1{\Omega_S^{#1}}

\def\otO{\otimes_{\cO}}
\def\cU{{\Cal U}}

\def\onab{\overline{\nabla}}

\def\HO{H_{\cO}}
\def\HQ{H_{\Bbb Q}}
\def\HC{H_{\Bbb C}}
\def\HQcU#1#2{\HQ^{#1}(\cU/S)(#2)}
\def\HCcU#1{\HC^{#1}(\cU/S)}
\def\HOcU#1{\HO^{#1}(\cU/S)}

\def\HcU#1#2{H^{#1,#2}(\cU/S)}

\def\HcU#1#2{H^{#1,#2}(\cU/S)}
\def\HcXp#1#2{H^{#1,#2}(\cX/S)_{prim}}

\def\HcUan#1#2{H^{#1,#2}(\cU/S)^{an}}

\def\HQcU#1#2{\HQ^{#1}(\cU/S)(#2)}

\def\HQcU#1#2{\HQ^{#1}(\cU/S)(#2)}

\def\HCcU#1{\HC^{#1}(\cU/S)}
\def\HCcX#1{\HC^{#1}(\cX/S)}
\def\HCcU#1{\HC^{#1}(\cU/S)}
\def\HCcXp#1{\HC^{#1}(\cX/S)_{prim}}

\def\HOcU#1{\HO^{#1}(\cU/S)}
\def\HOcX#1{\HO^{#1}(\cX/S)}
\def\HOcU#1{\HO^{#1}(\cU/S)}

\def\WSan#1{\Omega_{\San}^{#1}}

\def\Xx{X_x}
\def\Zx{Z_x}
\def\Ux{U_x}
\def\xb{\overline{x}}
\def\Uxb{U_{\xb}}

\def\Zxx{Z_{x}}
\def\Zst{Z}
\def\cZst{\cZ}

\def\Yst{Y}

\def\rocXZ{\kappa_{(\cX,\cZ)}}

\def\roolog{\kappa_0^{log}}

\def\roxlog{\kappa_x^{log}}

\def\TS{\Theta_S}
\def\TxS{T_x S}

\def\UC{U_{\Bbb C}}

\def\SC{S_{\Bbb C}}

\def\WS#1{\Omega_S^{#1}}
\def\WX#1{\Omega_X^{#1}}

\def\TX{T_X}

\def\TXZx{T_{\Xx}(-\log \Zxx)}

\def\WXZ#1{\Omega_X^{#1}(\log \Zst)}

\def\WcXZ#1{\Omega_{\cX/S}^{#1}(\log \cZst)}

\def\WcXkZ#1{\Omega_{\cX/k}^{#1}(\log \cZst)}

\def\WPnY#1{\Omega_{\Bbb P^n}^{#1}(\log \Yst)}
\def\TXZ{T_X(-\log \Zst)}
\def\TcXZS{T_{\cX/S}(-\log \cZst)}

\def\fXZ#1{\phi_{X,Z}^{#1}}
\def\cXZ{\psi_{(X,Z)}}

\def\cXZx{\psi_{(\Xx,\Zx)}}

\def\eXZ{\eta_{(X,Z)}}

\def\ccXZS{c_S(\cX,\cZ)}
\def\ccXZtS{c_{\wtd{S}}(\wtd{\cX},\wtd{\cZ})}

\def\tS{\wtd{S}}

\def\scs{\hbox{ }:\hbox{ }}
\def\onab{\overline{\nabla}}

\def\dlog#1{\frac{d#1}{#1}}
\def\chB#1#2#3{ch_{B,#3}^{#1,#2}}
\def\chD#1#2#3{ch_{D,#3}^{#1,#2}}
\def\chet#1#2#3{ch_{et,#3}^{#1,#2}}
\def\chcont#1#2#3{ch_{cont,#3}^{#1,#2}}
\def\Ql{\Bbb Q_{\ell}}
\def\Zl{\Bbb Z_{\ell}}
\def\etab{\overline{\eta}}

\def\Pd{\overset{\vee}\to{\Bbb P}}
\def\kb{\overline{k}}
\def\xb{\overline{x}}

\def\pitop{\pi_1^{top}}
\def\pialg{\pi_1^{alg}}

\def\wcU{\omega_{\cU/S}(\sigma)}
\def\dcU{\delta_{\cU/S}(\sigma)}

\def\ccU{c_{\cU/S}(\sigma)}
\def\cUx{c_{\Ux}(\sigma)}
\def\cUo{c_{U}(\sigma)}
\def\cUetab{c_{U_{\etab}}(\sigma)}

\def\Dqi{\Delta_{\gamma_i}}
\def\Dq{\Delta_{\gamma}}
\def\spa{\hbox{ }}

\def\nabmo{\overline{\nabla}^{m,0}}
\def\nab#1#2{\overline{\nabla}^{#1,#2}}

\def\Pol{P}
\def\h#1#2{h_{#1}(#2)}

\def\San{S_{an}}
\def\Szar{S_{zar}}

\def\dt{\partial_{tame}}
\def\dd{\partial_{div}}

\def\dz{\partial_Z}

\head \S5. Implication on injectivity of Chern class maps for $K_1$ of
surfaces. \endhead
\vskip 8pt

Let $X$ be a projective smooth surface over a field $k$.
Let $U\subset X$ be the complement of a simple normal crossing
divisor $Z\subset X$. In this section we discuss an implication of
the surjectivity of $\chB 22 U$ and $\chet 22 U$ on the Chern class maps
for $CH^2(X,1)$. Recall that $CH^2(X,1)$ is by definition the cohomology of
the following complex
$$ K_2(k(X)) @>{\dt}>> \bigoplus_{C\subset X} k(C)^*
@>{\dd}>> \bigoplus_{x\in X} \Bbb Z,$$
where the sum on the middle term ranges over all irreducible curves on $X$
and that on the right hand side over all closed points of $X$.
The map $\dt$ is the so-called tame symbol and $\dd$ is the sum of divisors
of rational functions on curves. Thus an element of $\Ker(\dd)$ is given by a
finite sum
$\sum_{i} (C_i,f_i)$, where $f_i$ is a non-zero rational function on
an irreducible curve $C_i\subset X$ such that $\sum_{i} div(f_i)=0$ on $X$.
We recall that $K_2(k(X))$ is generated as an abelian groups by symbols
$\{f,g\}$ for non-zero rational functions $f,g$ on $X$ and that
$$ \dt(\{f,g\})=((f)_0,g)+((f)_{\infty},1/g)+((g)_0,1/f)+((g)_{\infty},f),$$
where $(f)_0$ (resp. $(f)_{\infty}$) is the zero (resp. pole) divisor of $f$.
An important tool to study $CH^2(X,1)$ is the Chern class map:
In case $k=\Bbb C$ it is given by
$$ \chD 21 X \scs CH^2(X,1) \to H^3_D(X,\Bbb Z(2)),$$
where the group on the right hand side is the Deligne cohomology group
(cf. [EV] and [J2]).
Assuming the first Betti number $b_1(X)=0$, we have the following explicit
desciption of $\chD 21 X$. Take $\alpha=\sum_{i} (C_i,f_i)\in \Ker(\dd)$.
Under the isomorphism (cf. [EV, 2.10])
$$ H^3_D(X,\Bbb Z(2))\simeq \frac{H^2(X,\Bbb C)}{H^2(X,\Bbb Z(2))+ F^2H^2(X,\Bbb C)}
\simeq \frac{F^1H^2(X,\Bbb C)^*}{H_2(X,\Bbb Z)},$$
$\chD 21 X (\alpha)$ is identified with a linear function on complex valued
$C^{\infty}$-forms $\omega$ and we have
$$ \chD 21 X (\alpha)(\omega)=
\frac{1}{2\pi\sqrt{-1}} \sum_i \int_{C_i-\gamma_i} \log(f_i)\omega
+ \int_{\Gamma} \omega,$$
where $\gamma_i:=f_i^{-1}(\gamma_0)$ with $\gamma_0$, a path on
$\Bbb P^1_{\Bbb C}$ connecting $0$ with $\infty$ and $\Gamma$ is a real
piecewise smooth $2$-chain on $X$ such that
$\partial\Gamma =\cup_i \gamma_i$ which exists due to the assumption
$\alpha\in \Ker(\dd)$ and $b_1(X)=0$.
\vskip 5pt
In case $k$ is a finite extension of $\Bbb Q$ we have the chern class map
$$\chcont 21 X \scs CH^2(X,1)\otimes\Zl\to
H_{cont}^3(X,\Zl(2)),$$
where $H^i_{cont}$ denotes the continuous etale cohomology of $X$ (cf. [J3]).
\vskip 5pt

Now let $Z\subset X\supset U$ be as in the begining of this section and write
$Z=\cup_{1\leq i\leq r} Z_i$ with $Z_i$, smooth irreducible curves intersecting
transversally with each other.
We consider the following subgroup of $\Ker(\dd)$
$$ CH^1(Z,1)=\Ker(\bigoplus_{1\leq i\leq r} k(Z_i)^* @>{\dd}>>
\bigoplus_{x\in Z} \Bbb Z).$$
By the localization theory for higher Chow group we have the exact sequence
$$ CH^2(U,2) @>{\dz}>> CH^1(Z,1) \to CH^2(X,1)$$
where the first map coincides up to sign with the composite of
the natural map $CH^2(U,2) \to K_2(k(X))$ and $\dt$.

\Th 5-1. \it (1) Assume $k=\Bbb C$ and that there exists a subspace
$\Delta\subset CH^2(U,2)\otimes \Bbb Q$ such that the restriction of
$\chB 22 U$ on $\Delta$ is surjective. Let $\alpha\in CH^1(Z,1)$ and assume
$\chD 21 X(\alpha)=0$. Then $\alpha\in \dz(\Delta)$ in
$CH^1(Z,1)\otimes\Bbb Q$. In particular $\alpha=0$ in $CH^2(X,1)\otimes\Bbb Q$.
\vskip 4pt\noindent
(2) Assume that $k$ is a finite extension of $\Bbb Q$.
The analogous fact holds for $\chet 22 U$ and $\chcont 21 X$.
\rm\vskip 5pt

\Rem 5-1. \it
The main results in \S3 and \S4 imply that in case $X$ is a generic
hypersurface of degree $d\geq 4$ and $Z$ is the union of generic hypersurface sections on $X$, there exists $\Delta$ satisfying the assumption of Th.(5-1)
with such explicit generators as given in Def.(2-1)(2).
\rm\vskip 5pt

The following result is a direct consequence of Th.(5-1) and Th.(0-1).
Let $\cZ\subset \cX$ be as in the introduction and let $\Zx\subset \Xx$ be its
fibers over $x\in S$.

\Cor 5-2. \it Assume $\underset{1\leq i\leq r}\to{\sum} d_i \geq n+1+\ccXZS$
and $n-r=2$.
\roster
\item
Assume $k=\Bbb C$. There exists $E\subset S(\Bbb C)$ which is
the union of countable many proper analytic subset of $S(\Bbb C)$ such that
$\chD 21 {\Xx}$ restricted on the image of $CH^1(\Zx,1)$ is injective modulo
torsion for $\forall x\in S(\Bbb C)\setminus E$.
\item
Assume that $k$ is a finite extension of $\Bbb Q$ and
$S(k)\not=\emptyset$.
Let $\pi: S \to \Bbb P_k^N$ be a dominant quasi-finite morphism.
There exist a subset $H\subset\Bbb P_k^N(k)$ such that:
\vskip 4pt\noindent
$(i)$
$\chD 21 {\Xx}$ restricted on the image of $CH^1(\Zx,1)$ is injective modulo
torsion for any closed point $x\in S$ such that $\pi(x)\in H$.
\vskip 4pt\noindent
$(ii)$
Let $\Sigma$ be any finite set of primes of $k$ and let $k_v$ be the
completion of $k$ at $v\in \Sigma$. Then the image of $H$ in
$\prod_{v\in \Sigma} \Bbb P_k^N(k_v)$ is dense.
\endroster
\rm

\Rem 5-2. \it Let the notation be as in Cor.(5-2)(1).
In the forthcoming paper [AS2] it is shown that there exist
$x\in S(\Bbb C)\setminus E$ such that the image of $CH^1(\Zx,1)$ in
$CH^2(\Xx,1)$ is non-torsion. Thus Cor.(5-2) has indeed non-trivial implication
on the injectivity of the Chern class map.
\rm
\vskip 5pt\noindent
\it Proof of Th.(5-1). \rm
The idea of the following proof is taken from [J1, 9.8].
We only treat the first assertion. The second is proven by the same way.
We have the commutative diagram (cf. [Bl] and [J2, 3.3 and 1.15]).
$$\matrix
CH^2(X,2) & \to & CH^2(U,2) &\to& CH^1(Z,1) &\to& CH^2(X,1) \\
\downarrow\rlap{$\chD 22X$} && \downarrow\rlap{$\chD 22U$}&&
\downarrow\rlap{$\chD 11Z$} && \downarrow\rlap{$\chD 21X$}\\
H^2_D(X,\Bbb Z(2)) &@>\iota>>& H^2_D(U,\Bbb Z(2)) &\to&
H^3_{D,Z}(X,\Bbb Z(2)) &\to& H^3_D(X,\Bbb Z(2)) \\
\endmatrix$$
Here the horizontal sequences are the localization sequences
for higher Chow groups and Deligne cohomology groups and
they are exact. The vertical maps are Chern class maps.
By a simple diagram chasing it suffices to show that $\chD 11Z$ is
injective and
$\chD 22U(\Delta)+\Im(\iota)$ spans $H^2_D(U,\Bbb Q(2))$.
In order to show the first assertion we note the commutative diagram
$$\matrix
0\to&\underset{1\leq i\leq r}\to{\bigoplus} CH^1(Z_i,1) &\to& CH^1(Z,1) &\to&
\underset{x\in W}\to{\bigoplus} \Bbb Z\\
&\downarrow\rlap{$\chD 11{Z_i}$} && \downarrow\rlap{$\chD 11Z$}
&&  \llap{$\simeq$}\downarrow\rlap{$\chD 00x$}\\
0\to&\underset{1\leq i\leq r}\to{\bigoplus} H^1_{D}(Z_i,\Bbb Z(1)) &\to&
H^3_{D,Z}(X,\Bbb Z(2)) &\to&
\underset{x\in W}\to{\bigoplus} H^0_{D}(x,\Bbb Z(0))\\
&\downarrow\rlap{$\simeq$}&&\downarrow\rlap{$=$}&&\downarrow\rlap{$\simeq$}\\
&\underset{1\leq i\leq r}\to{\bigoplus} H^3_{D,Z_i}(X,\Bbb Z(2)) &\to&
 H^3_{D,Z}(X,\Bbb Z(2)) &\to&
\underset{x\in W}\to{\bigoplus}H^4_{D,x}(X,\Bbb Z(2))\\
\endmatrix$$
where $W=\coprod_{1\leq i\not=j\leq r} Z_i\cap Z_j$.
The horizontal sequences come from the Mayer-Vietoris spectral sequence
and they are exact. Thus the desired assertion follows from the fact that
$\chD 11{Z_i}$ is an isomorphism (cf. [J2, 3.2]).
To show the second assertion we recall the exact sequence
(cf. [EV, 2.10])
$$ 0 \to H^1(U,\Bbb C)/H^1(U,\Bbb Z(2)) @>\gamma>> H^2_D(U,\Bbb Z(2))
@>\beta>> H^2(U,\Bbb Z(2))\cap F^2H^2(U,\Bbb C) \to 0$$
and the same sequence with $U$ replaced by $X$. We have
$\chB 22 U=\beta\cdot \chD 22U$ and the commutative diagram
$$\matrix
\Bbb C^*\otimes CH^1(U,1) &\to& CH^1(U,1)\otimes CH^1(U,1)\\
\downarrow\rlap{$\log\otimes \chB 11 U$}&&\downarrow\rlap{\text{product}}\\
\Bbb C/\Bbb Z(1)\otimes H^1(U,\Bbb Z(1)) && CH^2(U,2)\\
\downarrow\rlap{$\simeq$}&&\downarrow\rlap{$\chD 22U$}\\
H^1(U,\Bbb C)/H^1(U,\Bbb Z(2)) &@>\gamma>>& H^2_D(U,\Bbb Z(2))\\
\endmatrix$$
Therefore the desired assertion follows from the surjectivity of
$$CH^1(U,1) @>{\chB 11 U}>> H^1(U,\Bbb Z(1))/H^1(X,\Bbb Z(1))$$
which is easily seen.
\qed

\vskip 20pt

\input amstex
\documentstyle{amsppt}
\hsize=16cm
\vsize=23cm

\document
\NoBlackBoxes
\nologo


\enddocument